\documentclass[12pt]{article}

\usepackage{amsmath,a4wide,amssymb,amsfonts,amsbsy,xy,latexsym,epsf,pstricks,times}
\xyoption{all}

\usepackage[latin1]{inputenc}

\usepackage{hyperref}

\usepackage{tikz}
\usetikzlibrary{matrix}
\usetikzlibrary{decorations}
\usetikzlibrary{decorations.pathreplacing}
\newenvironment{preuve}{\medbreak \noindent {\bf Proof~---}}
                       {\hfill $\square$ \medbreak}

\def\cB{X}
\def\cC{Y}
\def\cD{Z}
\def\cS{S}
\def\cE{E}
\def\cF{F}
\def\cJ{C}
\def\cK{D}
\def\cI{W}
\def\cG{G}

\def\cO{{\cal O}}
\def\cU{{\cal U}}

\def\PP{{\bf P}}
\def\AA{{\bf A}}
\def\NN{{\bf N}}
\def\ZZ{{\bf Z}}
\def\QQ{{\bf Q}}

\let\set\mathbb

\def\FF{{\set F}}
\def\SS{{\set S}}
\def\XX{{\set X}}
\def\YY{{\set Y}}
\def\MM{{\set M}}

\def\GG{{\set G}}

\def\Spec{\mathop{\rm{Spec}}\nolimits }
\def\Desc{\mathop{\rm{Desc}}\nolimits }
\def\Et{\mathop{\rm{Et}}\nolimits }
\def\Mod{{\set M}}
\def\Id{{\rm Id}}
\def\pgcd{\mathop{\rm{gcd}}\nolimits }
\def\Proj{\mathop{\rm{Proj}}\nolimits }

\def\Gal{\mathop{\rm{Gal}}\nolimits }
\def\Aut{\mathop{\rm{Aut}}\nolimits }

\def\btheta{{\bar \theta}}

\newtheorem{theoreme}{Theorem}[section]
\newtheorem{lem}[theoreme]{Lemma}
\newtheorem{definition}[theoreme]{Definition}
\newtheorem{proposition}[theoreme]{Proposition}

\newenvironment{remarque}{\begin{quote}
                     {{\bf Remark} --}
                    }{\end{quote}}

\begin{document}
\author{Jean-Marc Couveignes\thanks{This work is 
supported by the French Agence Nationale pour la Recherche through project ALGOL (ANR-07-BLAN-0248).} and  Emmanuel Hallouin\thanks{Institut de Math{\'e}matiques de Toulouse,
    Universit{\'e} de Toulouse et CNRS, Universit{\'e}  de Toulouse 2
le Mirail, 5 all{\'e}es
    Antonio Machado 31058 Toulouse c{\'e}dex 9}}
\title{Global   descent obstructions for  varieties}

\maketitle

\begin{abstract}
We show how to transport descent  obstructions  from the category
of covers to the category of varieties. We deduce examples of curves having $\QQ$ as field of moduli, that admit models
over every completion of $\QQ$, but have no model over $\QQ$.

AMS classification 11R34, 12G05, 11G35, 14D22.
\end{abstract}

\tableofcontents

\section{Introduction}\label{section:introduction}

If~$k$ is  a field, a $k$-variety is by definition
a separated 
scheme of finite type over $\Spec(k)$.
A $k$-curve is a variety of dimension $1$  over $k$. 
A $k$-surface is a variety of dimension $2$  over $k$. 

%C'est quoi un revêtement (courbes nodales, surfaces, ...)
%If $X$ and $Y$ are   $K$-varieties and $f : Y \rightarrow X$ a non-constant and separable morphism, we say the $f$ is a (possibly branched) covering of $X$.

\subsection{The main results}

In this article, we construct    descent obstructions   in the category of varieties.
For example, we show 
the following theorem:

\begin{theoreme}\label{theoreme:courbex}
There exists a  projective,  integral and smooth curve over~$\overline{\QQ}$,
having~$\QQ$ as field of moduli, which
is defined over all the completions of~$\QQ$ but not over~$\QQ$ itself.
\end{theoreme}

The main idea is to  start from  a 
 descent obstruction in
the category of covers of curves, and to transport it into various  other categories: the category of quasi-projective surfaces, the category of
proper surfaces, and finally the category of smooth curves.
This process is summarized by the following theorem:

\begin{theoreme}\label{theoreme:courbe}
Let $k$ be a field of characteristic zero, $k^a$ an algebraic closure of $k$.
Let~$X_k$ be a  smooth, projective, geometrically integral curve over~$k$ and let~$X$
denote the base change to~$k^a$ of~$X_k$. Let~$Y$ be a smooth, projective, integral curve
over~$k^a$ and let~$\varphi : Y \rightarrow X$ be a (possibly ramified) cover over~$k^a$,
having $k$ as field of moduli.
There exists a smooth, projective, integral curve over~$k^a$ having $k$ as field
of moduli and having exactly the same fields of definition as the initial cover  $\varphi$.
\end{theoreme}

Examples of obstructions to descent have been mostly
constructed in  the categories of $G$-covers and covers
\cite{CH, DF,  CG} and in the category
of dynamical systems \cite{sil}.
A key technical point is that, in many cases,
 one can measure these obstructions in terms of the  Galois 
cohomology of a  finite abelian
group.
As far as we know, no example of global  obstructions was
known for varieties. Mestre gave some examples
of local descent obstructions for hyperelliptic curves in~\cite{Mestre}.
Dèbes and Emsalem \cite{DE} give  a criterion for  a curve to be defined
over its field of moduli. This criterion  
involves a particular  model for the quotient of the curve by its
automorphism  group. Dèbes and Emsalem  prove that the local-global
principle applies to the descent problem for a  curve {\it together
with its automorphisms}. However they leave open the question of the
local-global principle for a curve (and a variety in general).

 Global  descent obstructions for covers
have been constructed
by Ros and  Couveignes:

\begin{theoreme}[cf.\cite{roscouveignes}, Corollaire~2] \label{theoreme:catrevetement}
There exists a connected ramified $\overline{\QQ}$-cover of~$\PP^1_\QQ$ having~$\QQ$ as field of moduli,
which is defined over all the completions of~$\QQ$ but which does not admit any model
over~$\QQ$.
\end{theoreme}

If we  apply theorem \ref{theoreme:courbe} to these obstructions, we 
prove  theorem~\ref{theoreme:courbex}.

\subsection{Overview of the paper}

Let $k$ be a field with characteristic zero.
Let $X_k$ be a smooth, projective and geometrically integral curve over $k$ and
set $X=X_k\times_{\Spec(k)}\Spec(k^a)$. The starting point of all the constructions, in the sequel,
is a smooth projective and integral curve~$Y$ over $k^a$ which covers~$X$, i.e. there
exists a non-constant morphism $\varphi : Y \rightarrow X$ of $k^a$-curves.
We would like to construct a variety having the same 
moduli and definition properties (the same
field of moduli  and the same fields of definition) as $\varphi$.
A first natural idea is to consider the complement $X\times Y-G(\varphi)$
of the graph $G(\varphi)$
of $\varphi$ in the product $X\times Y$. We call it the {\it mark} of $\varphi$.
We hope this surface will have the expected property: same field of moduli
and same fields of definition as $\varphi$. This will be true in many cases.
In order to prove it, we shall associate to $\varphi$ the stack of all
its models. We shall similarly associate to the mark of~$\varphi$ the stack of all
its models. We then try to  construct a morphism between these
two stacks. If this morphism happens to be an equivalence, then we are
done.

In section~\ref{section:gerbes}, we recall the definition of the stack and gerbe
of ``models'' of an algebraic variety over~$k^a$ (or of a cover of curves over~$k^a$). We then explain how a morphism between the two gerbes of models 
associated
with  two
objects relates the definition and module properties of either objects.
In the sequel we shall  make use of these functorial properties
to transport descent obstructions from a category to another one.
 It turns out that the key point is to control
the group of automorphisms of the involved objects.

To make this task easier, in section~\ref{section:sansauto}, we prove that we can suppose that
the base curve~$X$ of our starting cover~$\varphi$, does not
have any non trivial $k^a$-automorphism. In other words, we construct
another $k$-curve $X'_k$ without any non-trivial $k^a$-automorphism and a $k^a$-cover
$Y'\rightarrow X'_k\times_{\Spec(k)} \Spec(k^a)$ having the same field of moduli and the same
fields of definition as~$\varphi$.

In section~\ref{section:surfacesaffines}, we do suppose that~$X$ does not have any non-trivial
automorphism and we prove that the {\it mark} of $\varphi$ has the same field of moduli
and the same fields of definition as $\varphi$.

In section~\ref{section:surfacesnormalesprojectives}, we assume that the field of moduli of
the cover~$\varphi$ is~$k$ and we construct a  proper normal $k^a$-surface having $k$ as
field of moduli and the same fields of definition as $\varphi$. 
This proper surface is constructed as a cover
of $X \times Y$ which is strongly ramified along  the graph of $\varphi$.

Finally, in section~\ref{section:courbes}, we construct a projective $k^a$-curve, having~$k$
as field of moduli, and having the same fields of definition as the initial cover $\varphi$.
This curve is drawn on the previous  surface. It is  obtained by deformation of a stable
curve chosen to have the same automorphism group as the surface.

\medbreak

{\bf\noindent Notations.} If~$k$ is a field, we denote by~$k^a$ its algebraic closure.
Let~$l$ be a $k$-extension and let~$X_l$ be a $l$-variety. We denote
by~$\Aut_l(X_l)$ or simply~$\Aut(X_l)$, the group of automorphisms of the $l$-variety~$X_l$
(i.e. automorphisms over~$\Spec(l)$). On the other hand, we denote by~$\Aut_k(X_l)$
the group of automorphisms of the {\sl $k$-scheme}~$X_l$
(i.e. automorphisms over~$\Spec(k)$). For~$f \in l(X_l)$,
$(f)_0$ denote the divisor of zeros of the function~$f$ while~$(f)_\infty$
denote the divisor of poles of the function~$f$.

\section{Stack of ``models''}\label{section:gerbes}

In this section
$k$ is  a field of characteristic zero and 
$k^a$ is an algebraic closure of $k$.

\subsection{The conjugate of a variety}

Let $X$ be a $k^a$-variety.
We denote by $\pi : X\rightarrow  \Spec (k^a)$ the structural
morphism.
Let $\sigma : k^a\rightarrow k^a$ be a $k$-isomorphism.
We denote ${}^\sigma\! X$ the $k^a$-variety defined
to be $X$ itself with  the structural morphism
${}^\sigma\!\pi =    \Spec(\sigma)^{(-1)} \circ \pi$.
It is clear that the square below is cartesian.
So ${}^\sigma\!\pi$ is the pullback of $\pi$ along $\Spec(\sigma)$.

$$
\xymatrix@C=50pt@R=25pt{
{}^\sigma \! X \ar@{->}[r]^{\pi}\ar@{->}[d]^{\Id} & \Spec (k^a) \ar@{->}[r]^{\Spec(\sigma)^{(-1)}} 
& \Spec(k^a)\ar@{->}[d]^{\Spec(\sigma)}\\
 X \ar@{->}[rr]^{\pi} && \Spec (k^a) 
}
$$

With this (slightly abusive) notation one has ${}^\tau\!({}^\sigma\!(\pi ))
={}^{\tau\sigma}\! \pi$ and  ${}^\tau\!({}^\sigma\!(X ))
={}^{\tau\sigma}\! X$. If $X$ is an affine variety,  then ${}^\sigma\! X$
is obtained from $X$ by letting $\sigma$ act on the coefficients in the defining
equations of $X$. One may prefer to write
$X^{\Spec(\sigma)}$ rather than ${}^\sigma\! X$. This is fine also and we do have
$(X^{\Spec(\sigma)})^{\Spec(\tau)}=X^{\Spec(\sigma)\circ \Spec(\tau)}$.

\subsection{The field of moduli}

It is natural to ask if $X$ and ${}^\sigma\! X$ are isomorphic.
They are certainly isomorphic as schemes (and even equal by definition). But
 as varieties over $k^a$, they are isomorphic if and only if there exists
an isomorphism $\phi_\sigma$ that makes the following diagram  commute

$$
\xymatrix@=25pt{
X \ar@{->}[rrr]^{\phi_\sigma}\ar@{->}[ddr]_{\pi} &&& X \ar@{->}[dl]^{\pi} \\
 && \Spec(k^a) \ar@{->}[dl]^{\Spec(\sigma)^{(-1)}}& \\
 & \Spec(k^a) && 
}
$$

The above triangle gives rise to  a commutative square

\begin{equation}\label{eq:fieldmod}
\vcenter{
\xymatrix@C=35pt@R=25pt{
X \ar@{->}[r]^{\phi_\sigma}\ar@{->}[d]_{\pi}  & X \ar@{->}[d]^{\pi} \\
  \Spec(k^a) \ar@{->}[r]^{\Spec(\sigma)} &\Spec(k^a) 
}
}
\end{equation}

The existence of such a square means that the isomorphism $\Spec(\sigma)$
of $\Spec(k^a)$
lifts to an isomorphism $\phi_\sigma$ of $X$.
If there exists such a lift $\phi_\sigma$ for every $\sigma$ in the absolute
Galois group of $k$,  then we say that the condition {\it field of moduli} is met, or that
$X$ has $k$ as  field of moduli.

\subsection{Fields of definition}\label{sec:fielddef}

Another natural question is, given $l\subset k^a$  an algebraic extension of 
$k$, does there exist an $l$-variety $\pi_l : X_l\rightarrow \Spec(l)$ and a cartesian
square 

\begin{equation}\label{eq:fielddef}
\vcenter{
\xymatrix@C=35pt@R=25pt{
X_l \ar@{<-}[r]\ar@{->}[d]_{\pi_l}  & X \ar@{->}[d]^{\pi} \\
  \Spec(l) \ar@{<-}[r]^{\Spec(\subset)} &\Spec(k^a) 
}
}
\end{equation}
where the line downstairs is the spectrum of the inclusion.
If such a square exists we say that $l$ is a 
{\it field of definition} of $X$. We say that $\pi_l: X_l\rightarrow \Spec(l)$
is a {\it model} of $\pi : X \rightarrow \Spec(k^a)$ over $l$.
One may wonder if it is important to impose the arrow
downstairs in the definition above. The answer is yes
in general. The existence
of such a  cartesian square may depend on the chosen
arrow downstairs. However, if the condition {\it field 
of moduli} is met, then we may compose the cartesian squares in
\ref{eq:fieldmod} and \ref{eq:fielddef} 

\begin{equation*}%\label{eq:fielddef}
\xymatrix@C=35pt@R=25pt{
X_l \ar@{<-}[r]\ar@{->}[d]_{\pi_l}  & X \ar@{->}[d]^{\pi} & X \ar@{->}[l]_{\phi_\sigma}\ar@{->}[d]^{\pi} \\
  \Spec(l) \ar@{<-}[r]^{\Spec(\subset)} &\Spec(k^a) &\Spec(k^a)\ar@{->}[l]_{\Spec(\sigma)}
}
\end{equation*}
and choose the arrow downstairs we prefer.

Another simple observation: if $X$ has a model
$\pi_k : X_k\rightarrow \Spec(k)$ over $k$  then
the condition {\it field of moduli} is met. Indeed, we write
$X$ as a fiber product $X=X_k\times _{\Spec(k)}\Spec(k^a)$
and we  take for
$\phi_\sigma$ the fiber product $\Id_{X_k} \times_{\Spec(k)} \Spec(\sigma)$
where $\Id_{X_k} : X_k\rightarrow X_k$ is the identity on $X_k$.

One may ask if the converse is true.

\subsection{Descent obstructions}

Assume the condition {\it field of moduli} holds true. Does there exist
a model over $k$ ? If the answer is no, we say that there is a 
{\it descent obstruction}.
In case $k$ is a number field, we say that the obstruction is 
{\it global} if 
\begin{enumerate}
\item the condition {\it field of moduli} holds true,
\item there is no model over $k$,
\item  but for every place $v$
of $k$ there exists a model of $X$ over some   extension 
$l\subset k^a$ such that $l$ can be embedded in the completion
$k_v$ of $k$ at $v$. 
\end{enumerate}

\subsection{The fibered category of ``models'' of a variety}

We denote by $\Et / \Spec(k)$ the category of (finite) étale morphisms
over $\Spec(k)$. An object 
 $U$
in this category
is a structural morphism  $\Spec (l)\rightarrow \Spec (k)$
where $k \hookrightarrow l$ is a finite  {\'e}tale  $k$-algebra.
We define a {\it covering} of $U$ to be a surjective family
$(U_i\rightarrow U)_i$
of morphisms  in $\Et /\Spec(k)$. This turns $\Et /\Spec(k)$ into a site called
the étale site on $\Spec(k)$. It indeed satisfies the three axioms 
of site : the pullback of a covering exists and is a covering;
a covering of a covering is a covering; and the identity is a covering.

Note that in this paper, we use the word {\it covering} in the context
of sites. We keep the word {\it cover}  for a non-constant (separable)
morphism between two smooth projective and geometrically integral
curves.

Now given a $k^a$-variety $X$, we define the fibered category over $\Et / \Spec(k)$ of
its ``models''. So, for any $k$-algebra $l$ we must say  what we mean
by a ``model'' of $X$ over $\Spec(l)$.
We first assume that $l$ is a field. We say that an $l$-variety 
$\pi_l: X_l \rightarrow \Spec(l)$
is a ``model'' of $X$ over  $\Spec(l)$ if and only if there
exists an embedding $e  : l \hookrightarrow   k^a$ of $k$-algebras and
a cartesian square:
\begin{equation*}
\xymatrix@C=35pt@R=25pt{
X_l \ar@{<-}[r]\ar@{->}[d]_{\pi_l}  & X \ar@{->}[d]^{\pi} \\
  \Spec(l) \ar@{<-}[r]^{\Spec(e)} &\Spec(k^a) 
}
\end{equation*}
We insist that this time we do not fix an embedding
of $l$ into $k^a$. In particular, if $l$ is a subfield of $k^a$
containing $k$, we accept models of $X$ but also models of all its
conjugates. So the word model here is less restrictive than
in section \ref{sec:fielddef}. This is why 
 we write the word model
between quotation marks in that case. Of course, as already noticed,
the two notions are equivalent when the condition {\it field
of moduli} holds true.

If $l$ is any finite étale algebra
over $k$, then it is a direct product of finitely many
finite field extensions of $k$. We define a ``model'' of $X$ over $\Spec(l)$
to be a disjoint union of ``models'' of $X$ over every connected component
of $\Spec(l)$.

\begin{definition}[The category $\Mod_X$ of ``models'' of $X$] \label{def_Mod_X}
Let~$X$ be a $k^a$-variety. The category of ``models'' of $X$, denoted~$\Mod_X$,
is the category:
\begin{itemize}
\item whose objects are all ``models'' of $X$ over all finite étale $k$-algebras,
\item and whose morphisms are the cartesian squares
\begin{equation*}
\xymatrix@C=5pt@R=25pt{
X_l \ar@{->}[rr]\ar@{->}[d]_{\pi_l}  && X_m \ar@{->}[d]^{\pi_m} \\
  \Spec(l) \ar@{->}[rr]   \ar@{->}[dr] &&  \ar@{->}[dl]\Spec(m) \\
&\Spec(k)&
}
\end{equation*}
\end{itemize}
The functor that associates $\Spec(l)$  to every  ``model''
over $\Spec(l)$ turns $\Mod_X$ into a fibered category  over
$\Et/\Spec(k)$; we denote by~$\Mod_X(l)$ or~$\Mod_X(\Spec(l))$ the fiber over~$\Spec(l)$.
\end{definition}

In particular, we can pullback a ``model'' $X_l\rightarrow \Spec(l)$ along any morphism
$\Spec(m)\rightarrow \Spec(l)$ over $\Spec(k)$. Note that 
pulling back is not quite innocent since it can turn a model
into its conjugates so to say.

\subsection{Descent data} \label{s_descent}

In fact, under mild conditions, the fibered category $\Mod_X$
happens to be a stack. In order to see this, we need to recall 
a few definitions and elementary results about descent data
(see Giraud \cite{Gir} or the more accessible  notes \cite{Vis} by Vistoli).

Let $\SS$ be a site and let $\XX$ be a fibered category over
$\SS$. Let $U$ be an object in $\SS$ and let $\cU=
(U_i\rightarrow U)_i$
be a covering of $U$. A {\it descent datum} from
$\cU$ to $U$ is a collection of objects $X_i \to U_i$.
For every $i$ and every $j$ we also want
a morphism
$\phi_{ij} : \pi_2^*(X_j) \rightarrow \pi_1^*(X_i)$ where
$\pi_1$ and $\pi_2$ are the two ``projections''
 in the cartesian 
 diagram

\begin{equation*}
\xymatrix@=25pt{
&U_i\times_U U_j  \ar@{->}[dr]_{\pi_2} \ar@{->}[dl]^{\pi_1}& \\
U_i \ar@{->}[dr] &  &U_j  \ar@{->}[dl]\\
& U&}
\end{equation*}

We also require that the following compatibility relation holds true
for any $i$, $j$, and $k$

\begin{equation}\label{eq:compat}
\pi_{12}^*(\phi_{ij}) \circ \pi_{23}^*(\phi_{jk}) =
\pi_{13}^*(\phi_{ik})
\end{equation}
where the $\pi_{12}$, $\pi_{23}$, $\pi_{31}$ are the partial
``projections'' in the cube below

\begin{equation*}
\xymatrix@=25pt{
&U_{ijk} \ar@{->}[dl]_{\pi_{12}} \ar@{->}[rr]^{\pi_{23}}&  & U_{jk}\ar@{->}[dl] \ar@{->}[dd]\\
U_{ij}\ar@{->}[rr] \ar@{->}[dd]& & U_j\ar@{->}[dd] &\\
 & U_{ik}  \ar@{<-}[uu]^(.3){\pi_{13}}   \ar@{->}[dl]\ar@{->}[rr] & &U_k\ar@{->}[dl]\\
U_i \ar@{->}[rr] & & U &  
}
\end{equation*}
and $U_{ij}=U_i\times_UU_j$, $U_{ijk}=U_i\times_UU_j\times_UU_k$.

A morphism of descent data is a collection of local morphisms
that are compatible with the glueing morphisms on either sides.
We thus obtain a category $\Desc_\XX(\cU, U)$ for every
covering $\cU$ of $U$. We denote by $\XX (U)$ the fiber
of $\XX$ above $U$. There is a 
functor $\XX (U) \rightarrow \Desc_\XX(\cU, U)$ that 
associates to any object over $U$ the collection 
of its restrictions over the $U_i$.
These constructions are functorial. For example, if
$\YY$ is another fibered category  and $\FF : \XX \rightarrow
\YY$ a cartesian functor,  then $\FF$ induces
a functor from $\XX(U)$ to $\YY(U)$
and a functor from $\Desc_\XX(\cU, U)$ to $\Desc_\YY(\cU,U)$.
Further, 
the composite functors $\XX(U) \rightarrow \YY(U) \rightarrow
\Desc_\YY(\cU, U)$ and $\XX(U)\rightarrow \Desc_\XX (\cU,U)\rightarrow 
\Desc_\YY(\cU,U)$ are isomorphic.

A fibered category $\XX$ over a site $\SS$
is a {\it stack} if and only if all
the functors $\XX (U) \rightarrow \Desc_\XX(\cU, U)$ 
are equivalences of categories.

\subsection{When $\Mod_X$  is a stack, next a gerbe} \label{s_stack_gerbe}

If $X$ is a $k^a$-variety then $\Mod_X$ is a 
fibered category  over 
$\Et / \Spec(k)$ and it makes sense to ask if it is a stack.

We first notice that if
$l$ and $m$ are two finite field extensions of $k$  and
if $l\subset m$, then $\Spec(m)\rightarrow \Spec(l)$
is a covering of $\Spec(l)$. If further $m$  is
 a Galois extension of  $l$, then a descent datum
from $\Spec(m)$ to $\Spec(l)$ is a model $\pi_m : 
X_m\rightarrow \Spec(m)$ of $X$ over $\Spec(m)$ and, 
for every $\sigma $ in $\Gal(m/l)$, an automorphism
$\phi_\sigma : X_m\rightarrow X_m$ of {\sl $l$-scheme}, such that the following
diagram commutes
\begin{equation*}%\label{eq:fieldmod}
\xymatrix@C=35pt@R=25pt{
X_m \ar@{->}[r]^{\phi_\sigma}\ar@{->}[d]_{\pi_m}  & X_m \ar@{->}[d]^{\pi_m} \\
  \Spec(m) \ar@{->}[r]^{\Spec(\sigma)} &\Spec(m) 
}
\end{equation*}
We emphasize the fact that each $\phi_\sigma$ need not be an automorphism of the $m$-variety~$X_m$
but only an automorphism of the $l$-scheme~$X_m$. Let~$\Aut_l(X_m)$ denote the set of
automorphisms of the $l$-scheme~$X_m$. The compatibility condition (\ref{eq:compat}) states that
the map $\Spec(\sigma) \mapsto \phi_\sigma$ must be a
group homomorphism from $\Aut_{\Spec(l)}(\Spec(m))$ into
$\Aut_l(X_m)$.

\begin{proposition}
Let~$X$ be a variety over~$k^a$. 
If~$X$ is affine or projective or if every finite subset of $X(k^a)$ is contained in an affine
subvariety, then the fibered category~$\Mod_X$ is a stack over $\Et/\Spec(k)$.
\end{proposition}

\begin{preuve}
This is a consequence of Weil's descent theory. See the initial article of Weil \cite{Weil}
or Serre's book \cite[Chap~V,\S4]{Ser}.
\end{preuve}

Recall that a locally non-empty and locally connected stack is called a {\sl gerbe}. 
More precisely a stack $\XX$  over a site $\SS$ is a gerbe
if and only if
\begin{enumerate}
\item For every object $U$ in $\SS$ there exists a covering
$(U_i\rightarrow U)_i$ of $U$ such that the fibers over 
the $U_i$ are non-empty,
\item Given two objects $X \mapsto U$ and $Y\mapsto U$
above $U$ there exists a covering $(U_i\rightarrow U)_i$ 
such that for every $i$ the pullbacks $X\times_UU_i$
and $Y\times_UU_i$ are isomorphic over $U_i$,
\item For every object $U$ in $\SS$ the fiber $\XX(U)$
is a groupoid.
\end{enumerate}

The stack~$\Mod_X$ of "models'' of a variety~$X$ always satisfies conditions one and three,
while the second one holds true if and only if $k$ is the field of moduli of~$X$.

\subsection{The stack, next the gerbe, of ``models'' of a cover of a curves} \label{s_cover_of_curves}

Since the starting point of our construction is a cover of curves, we need to define the stack
of ``model'' of a cover of curves. So we adapt the notions presented in the preceding subsections to this
context.

Let $X_k$ be a smooth, projective, geometrically integral curve over $k$.
We set $X=X_k\times_{\Spec(k)}\Spec(k^a)$.
Let $Y$ be a smooth projective and integral curve over $k^a$ and
let $\varphi : Y \rightarrow X$ be a non-constant  morphism of $k^a$
varieties.
Since $k$ has zero characteristic, the morphism
$\varphi$ is separable. We say that $\varphi$ is a cover
of $X$. Note that we allow branch points.
An isomorphism between two covers
$\varphi : Y\rightarrow X$ and $\psi : Z\rightarrow X$ is an isomorphism
of $k^a$-varieties 
$ i : Y \rightarrow Z$ such that $\psi\circ i=\varphi$.

{\sl The conjugate of a cover ---} If $\sigma$ is a $k$-automorphism of $k^a$,
 the conjugate 
variety ${}^\sigma \! X$ is obtained from  $X$ by composing
the structural morphism on the left with $\Spec(\sigma)^{(-1)}$.
The same is true for $Y$. So any $k^a$-morphism $\varphi$ between
$X$ and $Y$ can be seen as a $k^a$ morphism ${}^\sigma \! \varphi$
 between   ${}^\sigma \! X$ and ${}^\sigma  Y$. Since $X$
is the fiber product of  $X_k$ and $\Spec(k^a)$ over $\Spec(k)$,
we have a canonical isomorphism $\phi_\sigma = \Id_{X_k}\times_{\Spec(k)}
\Spec(\sigma)$ between $X$ and ${}^\sigma \! X$. The composite map
$\phi_\sigma^{(-1)}\circ {}^\sigma \! \varphi$ is a  morphism
of $k^a$-varieties between ${}^\sigma  Y$  and $X$. We call
it the conjugate of $\varphi$ by $\sigma$. We may denote it
${}^\sigma \! \varphi$ also by abuse of notation. 

{\sl The field of moduli ---} We say that the condition {\it field of moduli} holds true for $\varphi$,
or that~$k$ is the field of moduli of~$\varphi$, if ${}^\sigma \! \varphi$ is isomorphic to $\varphi$
for every $k$-automorphism $\sigma$ of $k^a$.

{\sl Fields of definition, models ---} If $l\subset k^a$ is an algebraic extension of $k$ we set
$X_l=X_k\times_{\Spec(k)}\Spec(l)$. Let $Y_l$ be a smooth
projective and geometrically connected $l$-curve. Let 
$\varphi_l : Y_l\rightarrow
X_l$ be a non-constant (separable)  map.
If we lift $\varphi_l$ along the spectrum of the inclusion $l\subset k^a$
we obtain a morphism from $Y_l\times_{\Spec(l)}\Spec(k^a)$
onto $X=X_l\times_{\Spec(l)}\Spec(k^a)$. If this cover is isomorphic
to $\varphi : Y\rightarrow X$ then we say that $\varphi_l$ is a model
of $\varphi$ over $l$.

So it makes sense to ask if there exist (global) obstructions
to descent for covers of curves. It it proven in
\cite{roscouveignes} that this is indeed the case.

{\sl The fibered category  of ``models'' of a cover ---}
Given a finite {\'e}tale $k$-algebra $l$
we explain what we mean by a  ``model'' of $\varphi$ over $l$.

Assume first that $l$ is a finite field extension of $k$.
Set $X_l=X_k\times_{\Spec(k)}\Spec(l)$. Let $Y_l$ be 
a smooth projective and geometrically integral curve over $\Spec(l)$
and $\varphi_l : Y_l\rightarrow X_l$ be a non-constant morphism
over $\Spec(l)$. 
We pick any
embedding $e : l\rightarrow
k^a$ of $k$-algebras. The pullback of $X_l$ along
$\Spec(e)$ is  $X$ (up to unique isomorphism) and we have the following
diagram

\begin{equation*}%\label{eq:fieldmod}
\xymatrix@=30pt{
&&Y_l\ar@{->}[dl]_{\varphi_l} \ar@{->}[ddl]&&Y_l\times_{\Spec(l)}\Spec(k^a)\ar@{->}[ddl]\ar@{->}[dl]\ar@{->}[ll]\\
X_k \ar@{->}[d]&X_l\ar@{->}[l] \ar@{->}[d]& &X\ar@{->}[ll]\ar@{->}[d]&\\
\Spec(k)&\Spec(l)\ar@{->}[l]&&\Spec(k^a)\ar@{->}[ll]_{\Spec(e)} &}
\end{equation*}

We say that $\varphi_l$ is a ``model'' of $\varphi$ if the cover
$$\varphi_l\times_{\Spec(l)}\Spec(k^a) : Y_l\times_{\Spec(l)}\Spec(k^a)
\rightarrow X$$ is isomorphic to $\varphi$.

Again we don't care about the choice of the embedding $e$. We just 
ask that such an embedding exists.

If $l$ is any finite {\'e}tale algebra
over $k$, we define a ``model'' of $\varphi$ over $\Spec(l)$
to be a disjoint union of ``models'' of $\varphi$ over every connected component
of $\Spec(l)$. We write $\Mod_\varphi$ for the category of all models of
$\varphi$. This is a fibered category over $\Et/\Spec(k)$.

The following proposition is a consequence
of  Weil's descent theorem:

\begin{proposition}
Let $X_k$ be a smooth, projective, geometrically integral curve over $k$ and
set $X=X_k\times_{\Spec(k)}\Spec(k^a)$.
Let $Y$ be a smooth projective and integral curve over $k^a$ and
let $\varphi  : Y \rightarrow X$ be a non-constant  morphism of $k^a$ curves.
Then the fibered category~$\Mod_\varphi$ is a stack over $\Et/\Spec(k)$.
\end{proposition}

Like for varieties, the stack~$\Mod_f$ is a gerbe if and only if $k$ is the field of moduli of~$f$.

\subsection{Transporting obstructions}

In this subsection by {\sl cover of curves}, we mean a cover~$f : Y \to X$
satisfying the hypotheses of section~\ref{s_cover_of_curves}.

Let us emphasize the following easy facts from the preceding sections.

\begin{proposition}
Let~$X$ be a $k^a$-variety (or a cover a curves) then:
\begin{enumerate}
\item the field~$k$ is the field of moduli of~$X$ if and only if the stack~$\Mod_X$ is a Gerbe;
\item the field~$l$ is a field of definition of~$X$ if and only if the fiber~$\Mod_X(l)$
is not empty.
\end{enumerate}
\end{proposition}

Let~$X$ and~$Y$ be two $k^a$-varieties. Using a (cartesian)
morphism of stacks from~$\Mod_X$ to~$\Mod_Y$,
we are now able to transport descent obstruction for~$X$ to descent obstruction for~$Y$.
Recall a cartesian morphism is a functor of fibered categories
that transforms cartesian square into cartesian squares.
So a cartesian  morphism~$\FF : \Mod_X \to \Mod_Y$,  associates
 an $l$-model~$\FF(X_l)$
of~$Y$ to every $l$-model~$X_l$ of~$X$, and commutes with base change.

\begin{proposition} \label{obstruction_transportor}
Let~$X$ and~$Y$ be either $k^a$-varieties or covers of curves. Suppose that there exists
a morphism~$\FF : \Mod_X \longrightarrow \Mod_Y$ of stacks.
\begin{enumerate}
\item If~$X$ has $k$ as field of moduli then $Y$ has~$k$ as field of moduli;
\item If~$l$ is a field of definition of~$X$ then~$l$ is also a field of definition of~$Y$.
\end{enumerate}
Moreover, if the first condition holds true and if~$\FF$ is fully faithful then:
\begin{enumerate}\setcounter{enumi}{2}
\item $l$ is a field of definition of~$X$ if and only if~$l$ is a field of definition of~$Y$.
\end{enumerate}
In that case, there is a descent obstruction for~$X$ if and only if there is one for~$Y$.
\end{proposition}

\begin{preuve}
The two first conditions are easy consequences of the preceding result. The third one can
be deduced form the following more general lemma.
\end{preuve}

\begin{lem}
Let~$\XX$ and~$\YY$ be two gerbes over a site~$\SS$ and let~$\FF : \XX \to \YY$ be a cartesian morphism.
If~$\FF$ is fully faithful then $\FF$ is essentially surjective.
\end{lem}

\begin{preuve}
Let~$U$ be an object in~$\SS$ and~$Y \to U$ an object in the fiber~$\YY(U)$. Locally~$\XX(U)$
is not empty: there exists a covering~$(U_i \to U)_i$ of~$U$ and objects~$X_i \in \XX(U_i)$
for all~$i$. Set~$Y_i = X \times_U U_i$. Locally,~$Y_i$ and~$\FF(X_i)$ are isomorphic: there
exists a covering~$(U_{ij} \to U_i)_j$ such that~$Y_i \times_{U_i} U_{ij}$
and~$\FF(X_i \times_{U_i} U_{ij})$ are isomorphic.
Set~$X_{ij} = X_i \times_{U_i} U_{ij}$ and~$Y_{ij} = Y_i \times_{U_i} U_{ij}$.

Note that the collection of objects~$(Y_{ij} \to U_{ij})_{ij}$ defines a descent datum
from~$(U_{ij} \to U)_{ij}$ to~$U$; indeed for every~$i,j, i',j'$, pulling back
identity gives rise to isomorphisms:
$$
\Phi_{iji'j'} : Y_{ij} \times_{U_{ij}} U_{iji'j'} \longrightarrow Y_{i'j'} 
\times_{U_{i'j'}} U_{iji'j'} 
$$
which clearly satisfy the compatibility conditions~(\ref{eq:compat}) of~\S\ref{s_descent}.

Since~$\FF$ is fully faithful, there exist isomorphisms
$$
\Psi_{iji'j'} :  X_{ij} \times_{U_{ij}} U_{iji'j'} \longrightarrow X_{i'j'} 
\times_{U_{i'j'}} U_{iji'j'} 
$$
which, in turn, satisfy the compatibility conditions~(\ref{eq:compat}) of~\S\ref{s_descent}.
We deduce that there exists~$X \to U$ in~$\XX(U)$ such that~$\FF(X) = Y$.
\end{preuve}

We end this section by an example of morphism between two stacks of ``models'' of a
variety.

%\begin{lem}\label{lem:quotient_by_subgroup}
%Let~$X$ be a $k^a$ variety having~$k$ as field of moduli.
%Suppose that the group~$\Aut(X)$ if finite and let~$G$ be a characteristic
%subgroup of~$\Aut(X)$. Then there exists  a cartesian
morphism from~$\Mod_X$ to~$\Mod_{X/G}$, where~$X/G$
%denotes the quotient variety of~$X$ by~$G$.
%\end{lem}
%
%\begin{preuve}
%We first need to define the image of an object. Let~$l$ be an extension of~$k$ and let~$X_l$
%be a ``model'' of~$X$ in~$\Mod_X(l)$. All the $k^a$-automorphisms
%of~$X_l \times_{\Spec(l)} \Spec(k^a)$ may not be defined over~$l$, but there exists a
%finite Galois extension~$m$ of~$l$ over which they are. Put~$X_m = X_l \times_{\Spec(l)} \Spec(m)$.
%One can now consider the quotient of~$X_m$ by the group~$G$: let~$p_m : X_m \to X_m/G$ be the
%canonical projection.
%
%Of course~$X_l$ is a model of~$X_m$ over~$l$; by section~\ref{s_stack_gerbe},
%there exists a group homomorphism~$\sigma \mapsto \phi_\sigma$ from~$\Gal(m/l)$
%to~$\Aut_l(X_m)$.
%
%For every~$g \in \Aut(X_m)$, and every~$\sigma \in \Gal(m/l)$, one
%has~$\phi_\sigma \circ g \circ \phi_\sigma^{-1} \in \Aut(X_m)$. Since~$G$ is a characteristic
%subgroup of~$\Aut(X_m)$, one also has~$\phi_\sigma \circ g \circ \phi_\sigma^{-1} \in G$,
%for every~$g \in G$. We deduce that~$p_m \circ \phi_\sigma \circ g = p_m \circ \phi_\sigma$
%for every~$g \in G$ so~$\phi_\sigma$ factorizes into~$\psi_\sigma : X_m/G \to X_m/G$.
%By uniqueness of this factorization, the correspondence~$\sigma \mapsto \psi_\sigma$
%is necessarily a group homomorphism from~$\Gal(m/l)$ to~$\Aut_l(X_m/G)$; therefore
%the quotient~$X_m/G$ descents to~$l$.
%
%Next, we need to define the image of a morphism [... A FINIR...]
%\end{preuve}

\begin{proposition}\label{prop:quotient_by_subgroup}
Let~$X$ be a integral variety over~$k^a$ having~$k$ as field of moduli and let~$G$ be a
finite subgroup of~$\Aut_{k^a}(X)$ which is normal in the group~$\Aut_k(X)$.
Assume that every orbit of $G$ is contained in an affine open subset
of $X$.
Then there is a morphism from~$\Mod_X$ to~$\Mod_{X/G}$, where~$X/G$
denotes the quotient variety of~$X$ by~$G$.

Let $Y \subset X/G$ be the complement of the branch locus of~$X \to X/G$. This an open
sub-variety of~$X/G$ and there is a morphism of stacks from~$\Mod_X$ to~$\Mod_Y$.
\end{proposition}

\begin{preuve}
We first need to define the image of an object. Let~$l$ be an extension of~$k$ and let~$X_l$
be a ``model'' of~$X$ in~$\Mod_X(l)$. All the elements of~$G$ may not be defined over~$l$,
but there exists a
finite Galois extension~$m$ of~$l$ over which they are. Put~$X_m = X_l \times_{\Spec(l)} \Spec(m)$.
One can now consider the quotient of~$X_m$ by the group~$G$: let~$p_m : X_m \to X_m/G$ be the
canonical projection. 
%Note that such  a quotient exists and is a variety 

Of course~$X_l$ is a model of~$X_m$ over~$l$; by section~\ref{s_stack_gerbe},
there exists a group homomorphism~$\sigma \mapsto \phi_\sigma$ from~$\Gal(m/l)$
to~$\Aut_l(X_m)$.
%%\begin{equation*}%\label{eq:fieldmod}
%%\xymatrix@=25pt{
%%X_m \ar@{->}[r]^{\phi_\sigma}\ar@{->}[d]_{\pi_m}  & X_m \ar@{->}[d]^{\pi_m} \\
%%  \Spec(m) \ar@{->}[r]^{\Spec(\sigma)} &\Spec(m) 
%%}
%%\end{equation*}
%%Moreover the~$\phi_\sigma$'s satisfy the usual compatibility conditions.

Since~$G$ is a normal subgroup of~$\Aut_k(X)$, it is, a fortiori, a normal subgroup
of~$\Aut_l(X_m)$
and thus, for every~$g \in G$ and
every~$\sigma \in \Gal(m/l)$, one has~$\phi_\sigma \circ g \circ \phi_\sigma^{-1} \in G$.
We deduce that~$p_m \circ \phi_\sigma \circ g = p_m \circ \phi_\sigma$
for every~$g \in G$. This implies that~$\phi_\sigma$ factorizes
into~$\psi_\sigma : X_m/G \to X_m/G$.
By uniqueness of this factorization, the correspondence~$\sigma \mapsto \psi_\sigma$
is necessarily a group homomorphism from~$\Gal(m/l)$ to~$\Aut_l(X_m/G)$; therefore
the quotient~$X_m/G$ descents to~$l$.

Next, we need to define the image of a morphism. Let~$X_i \to \Spec(l_i)$, $i = 1,2$ two
``models'' of~$X$. One can complete a cartesian square involving the~$X_i$'s in the
following way:
\begin{equation*}
\xymatrix@=25pt{
X_1 \ar@{->}[d] & X_2 \ar@{->}[d] \ar@{->}[l] & X_m \ar@{->}[d] \ar@{->}[l] \\
\Spec(l_1) & \Spec(l_2) \ar@{->}[l] & \Spec(m) \ar@{->}[l] \\
}
\end{equation*}
where $m$ is a finite Galois extension of $k$  such that all
 elements of~$G$ are defined over~$m$. We know that
there exist isomorphisms~$\Phi_1, \Phi_2$ making the following diagrams
commute:
$$
\xymatrix@=25pt{
Y_1 \ar@{->}[d] & X_m/G \ar@{->}[d] \ar@{->}[l]_{\Phi_1} \\
\Spec(l_1) & \Spec(m) \ar@{->}[l] \\
}
\qquad
\xymatrix@=25pt{
Y_2 \ar@{->}[d] & X_m/G \ar@{->}[d] \ar@{->}[l]_{\Phi_2} \\
\Spec(l_2) & \Spec(m) \ar@{->}[l] \\
}
$$
The image of the starting cartesian square is nothing else than:
$$
\xymatrix@=25pt{
Y_1 \ar@{->}[d] & Y_2 \ar@{->}[d] \ar@{->}[l]_{\Phi_1 \circ \Phi_2^{-1}} \\
\Spec(l_1) & \Spec(l_2) \ar@{->}[l] \\
}
$$
This completes the proof of the first statement.

The second statement is true because taking the branch locus commutes with base changes.
\end{preuve}

%In fact it is not necessary to have an equivalence of categories between~$\Mod_X$ and~$\Mod_Y$
%to share the fields of definition.

\section{Cancellation of the automorphism group of the base curve}\label{section:sansauto}

In this section, $k$ is a field of characteristic zero,~$k^a$ an algebraic closure of
it,~$l \subset k^a$ an algebraic extension of~$k$.
Let $X_k$ be a projective, smooth, geometrically integral curve over $k$
and set $X=X_k\times_{\Spec(k)}\Spec(k^a)$. We assume we are
given a smooth projective and integral curve~$Y$ over $k^a$ and a
cover~$\varphi : Y \to X$ having~$k$ as field of moduli.

We want to construct other covers having the same field of moduli
and the same fields of definition as~$\varphi$ but satisfying additional properties.
In particular, we want to show that one can assume that the base curve~$X$ has no
non-trivial $k^a$-automorphism.

We first  prove that the degree of the cover can be multiplied by any  prime integer not dividing 
the
initial degree.

\begin{proposition} \label{prop:granddegre}
Let $X_k$ be a  smooth, projective, geometrically integral curve over $k$
and set $X=X_k\times_{\Spec(k)}\Spec(k^a)$. Let~$Y$ be a smooth projective and integral
curve over $k^a$ and let~$\varphi : Y \to X$ be a cover over~$k^a$
of degree~$d$.
For every  prime~$p$ not dividing~$d$, there exists a smooth projective curve~$Y'$ over~$k^a$ and
a cover~$\psi : Y' \to X$ of degree~$pd$,  having the same field of moduli
and the same fields of definition as~$\varphi$.
\end{proposition}

\begin{preuve}
Let~$f \in k(X_k)$ be a non-constant function whose divisor is simple and does not meet
the ramification locus of~$\varphi$. The equation~$h^p = f$ defines a degree~$p$
extension of~$k(X_k)$. We denote by~$X'_k$ the smooth, projective, geometrically integral curve
corresponding to this function field
and we set~$X' = X'_k \times_{\Spec(k)} \Spec(k^a)$. The
morphism~$\nu : X' \to X$ is a cyclic Galois cover of degree~$p$.
We fix an algebraic closure~$\Omega$ of~$k^a(X)$ and embeddings of~$k^a(X')$ and~$k^a(Y)$
in~$\Omega$. 
Let~$Y'$ be the smooth projective $k^a$-curve corresponding to the
compositum of~$k^a(Y)$ and~$k^a(X')$. Since the field extensions~$k^a(Y)$ and~$k^a(X')$
are linearly disjoint over~$k^a(X)$, the cover~$\psi : Y' \to X$ has degree~$pd$:
$$
\xymatrix@=25pt{
& Y' \ar@{->}[dl]_{\varphi'} \ar@{->}[dr]^{\nu'} \ar@{->}[dd]^{\psi} \\
X' \ar@{->}[dr]_{\nu} & & Y \ar@{->}[dl]^{\varphi} \\
& X
}
$$
Let us prove that this construction yields a morphism of
stacks~$\FF : \Mod_\varphi \to \Mod_\psi$.

Let~$l \subset k^a$ be a finite extension of~$k$.
Set~$X_l = X_k \times_{\Spec(k)} \Spec(l)$, $X'_l = X'_k \times_{\Spec(k)} \Spec(l)$ and
consider~$\varphi_l : Y_l \to X_l$ an $l$-model of~$\varphi$.
In the construction above, one can replace~$X$, $X'$, $Y$
by~$X_l$, $X'_l$, $Y_l$. The $l$-curve~$Y'_l$ corresponding to the compositum of the two
function fields~$l(X'_l)$ and~$l(Y_l)$ is smooth, projective, geometrically integral
(because~$l$ is algebraically closed in the compositum) curve
and the $l$-cover~$\psi_l : Y'_l \to X_l$
is an $l$-model of~$\psi$. We define the morphism~$\FF$ by putting~$\FF(\varphi_l) = \psi_l$.
Since the function~$f$ has been chosen in~$k(\cB)$, the functor~$\FF$ maps cartesian squares to
cartesian squares. Thus~$\FF$ is a morphism of stacks. By
proposition~\ref{obstruction_transportor}, if~$l$ is a field
of definition of~$\varphi$ then $l$ is a field
of definition of~$\psi$ and if~$\varphi$ has $k$ as field of moduli then $\psi$
has $k$ as field of moduli.

To prove the converse, we use proposition~\ref{prop:quotient_by_subgroup} in order to construct a
morphism the other way around. Let~$\nu'$ denote the Galois cover~$Y' \to Y$.  We need
to show that the group~$\Aut(\nu')$ is normal in~$\Aut_k(\psi)$.
Let~$\Phi' \in \Aut_k(\psi)$. It induces maps~$\Phi : Y \to Y$ and~$\Psi : X \to X$
such that the following diagram commute:
\begin{equation} \label{diagram:genredegre}
\vcenter{
\xymatrix@C=5pt@R=15pt{
\cC' \ar@{->}[rr]^{\Phi'} \ar@{->}[d]_{\nu'} \ar@/_30pt/[dd]_{\psi} & & \cC'  \ar@{->}[d]^{\nu'} \ar@/^30pt/[dd]^{\psi} \\
\ar@{->}[d]_{\varphi} \cC \ar@{->}[rr]^\Phi&   & \ar@{->}[d]^{\varphi} \cC   \\
\cB\ar@{->}[rr]^{\Psi} \ar@{->}[dr]& & \cB \ar@{->}[dl] \\
& \Spec(k) &
}
}
\end{equation}
(horizontal arrows are morphisms of $k$-schemes). The existence of~$\Psi$ is a
consequence of the fact that~$X$ is defined over~$k$. The morphism~$\Phi$ exists
because~$Y \overset{\varphi}{\to} X$ is the maximal
sub-cover of~$Y' \overset{\psi}{\to} X$
unramified at the support of~$f$. And $f$ is $k$-rational.
Now if~$\Lambda \in \Aut_{k^a}(\nu')$,
i.e.~$\nu' \circ \Lambda = \nu'$, then:
$$
\nu' \circ \Phi' \circ \Lambda = \Phi \circ \nu' \circ \Lambda
= \Phi \circ \nu'= \nu' \circ \Phi'
$$
so~$\Phi' \circ \Lambda \circ \Phi'^{-1} \in \Aut_{k^a}(\nu')$, which was to be proven.
In conclusion, we do have a morphism~$\GG : \Mod_\psi \to \Mod_\varphi$ of stacks
and the lemma follows.
\end{preuve}

\begin{remarque}
The functor~$\FF : \MM_\varphi \to \MM_\psi$ is not fully faithful because~$\psi$ has more
automorphisms than~$\varphi$.
This is why, we do not apply point~$(3)$ in proposition~\ref{obstruction_transportor} here.
We instead construct another functor~$\GG : \MM_\psi \to \MM_\varphi$ and we apply points~$(1)$
and~$(2)$ in proposition~\ref{obstruction_transportor} to either functors~$\FF$ and~$\GG$
successively. We notice that~$\GG$ is a left inverse of~$\FF$.
\end{remarque}

Next, we show that  the base curve can be assumed
to have genus  greater than~$2$.

\begin{proposition} \label{prop:grandgenre}
Let $X_k$ be a  smooth, projective, geometrically integral curve over $k$
and set $X=X_k\times_{\Spec(k)}\Spec(k^a)$. Let~$Y$ be a smooth projective and integral
curve over $k^a$ and let~$\varphi : Y \to X$ be a cover over~$k^a$
of degree~$d$.
There exists a smooth, projective, geometrically integral curve~$\cB'_k$ over~$k$ of
genus greater than~$2$ and a cover~$\varphi' : Y' \to X'_k \times_{\Spec(k)} \Spec(k^a)$
having the same field of moduli and the same fields of definition as~$\varphi$.
\end{proposition}

\begin{preuve}
We use the construction and notation of diagram~\ref{diagram:genredegre} above.
We further assume that the chosen function has degree at least~$f$
has degree at least~$3$. By Hurwitz genus formula, the curve~$X'$ has a genus greater than or
equal to~$2$.

This construction yields a morphism of stacks~$\FF : \Mod_\varphi \to  \Mod_{\varphi'}$.
The cover~$\varphi : Y \to X$ is the maximal sub-cover of~$\psi : Y' \to X$ unramified at
the support of~$f$. Therefore, there exists a morphism
from~$\Aut_{k^a}(\varphi') \to \Aut_{k^a}(\varphi)$. This morphism is bijective
because~$k^a(X')$ and~$k^a(Y)$ are linearly disjoint over~$k^a(X)$. So the morphism~$\FF$
is fully faithful. We conclude, this time, using proposition~\ref{obstruction_transportor}.
\end{preuve}

Last, we prove that one can assume
the base curve to have no non-trivial $k^a$-automorphism.

\begin{proposition} \label{prop:sansauto}
Let $X_k$ be a  smooth, projective, geometrically integral curve over $k$
and set $X=X_k\times_{\Spec(k)}\Spec(k^a)$. Let~$Y$ be a smooth projective and integral
curve over $k^a$ and let~$\varphi : Y \to X$ be a cover over~$k^a$.
There exists a smooth, projective, geometrically integral curve~$\cB'_k$ over~$k$, of genus
greater that~$2$, such that~$\cB' = \cB'_k \times_{\Spec(k)} \Spec(k^a)$ does not have
any non-trivial automorphism and there exists a cover~$\varphi' : Y' \to X'$ over~$k^a$
having the same field of moduli and the same fields of definition as~$\varphi$.
\end{proposition}

\begin{preuve}
Thanks to proposition~\ref{prop:grandgenre}, one can assume
 that the genus~$g(X)$ of~$\cB$ is greater
then~$2$. Consequently, the group~$\Aut(\cB)$ of $k^a$-automorphisms is finite.

Let $p\ge 3$ be a prime integer. To begin with, we show that there exists a
non-constant function~$f \in k(\cB)$ which is non-singular above~$2$, $-2$ and $\infty$,
of degree greater than~$2+4p(g(\cB)-1)+2p^2$, such that the set~$f^{-1}(\{-2,2\})$
is not invariant by any non-trivial automorphism of~$\cB$,
and such that the set of singular values of~$\varphi$ does not meet the
set~$f^{-1}(\{2,-2,\infty \})$.

Let~$D$ be a simple effective divisor on~$X$ with degree greater than~$2+4p(g(\cB)-1)+2p^2$.
We also assume that~$D$ is disjoint from the set of singular values of~$\varphi$
and the linear space~$L(D)$ associated with~$D$ generates~$k^a(X)$ over~$k^a$.
In particular, for every~$\theta \in \Aut(X)$, this linear space is not
contained in the Kernel of~$\theta - \Id$. It is not contained in the kernel
of~$\theta + \Id$ either because it contains~$k^a$. If~$D$ has been chosen with a large
enough degree, the functions in~$L(D)$ having degree
less than the degree of~$D$ are contained in a finite union of strict
vector subspaces. Therefore there exists a non constant function~$f \in L(D)$
such that~$\deg(f) = \deg(D)$ and~$\theta(f) \not= \pm f$ for
all~$\theta \in \Aut(\cB) \setminus \{\Id\}$. By construction, this function is not
singular above~$\infty$ and~$f^{-1}(\infty)$ does not meet the singular values of~$\varphi$.
We can also assume that~$f \in k(X_k)$.

By construction, the function~$f^2$ does not have any non-trivial automorphism
(in short one has $\Aut_{k^a(f^2)}(k^a(\cB)) = \{\Id\}$). Using lemma~\ref{lem:autofibres},
we deduce that almost all the fibers of~$f^2$ are non singular and not fixed by any
non trivial automorphism of~$\Aut(\cB)$. In particular, there exists~$\lambda \in k^*$ such that
the fiber of~$f^2$ above~$\lambda^2$ is non singular, not fixed by any non-trivial
automorphism in~$\Aut(\cB)$ and does not meet the singular values of~$\varphi$.
The function~$2f/\lambda$ satisfies all the properties we want. Let us denote it by~$f$.

Now the equation~$h^p+h^{-p}-f=0$ defines a regular extension of~$k(\cB_k)$.
Let~$\cB''_k$ be the smooth, projective, geometrically integral curve associated to this function
field.
We denote by~$w$ the automorphism of~$\cB''_k$ given by~$w(h)=h^{-1}$
and by~$\cB'_k$ the quotient~$\cB''_k / \langle w\rangle$; this is a smooth, projective,
geometrically integral $k$-curve, covering~$\cB_k$ by a
$k$-cover~$\nu_k : \cB'_k\rightarrow \cB_k$ of degree~$p$.
Extending scalars to~$k^a$, we obtain a Galois cover~$\cB'' \to \cB$ of~$k^a$-curves,
with Galois group~$D_p$, and whose singular values are exactly $f^{-1}(\{2,-2,\infty \})$.
Since the subgroup~$\langle w \rangle$ is self-normalized in~$D_p$,
the quotient by this subgroup is a sub-cover~$\nu : \cB' \to \cB$ of $k^a$-curves
of degree~$p$ which does not have any non-trivial automorphism.

Because the ramification loci do not meet, the function fields~$k^a(\cB'')$ and~$k^a(\cC)$
are linear disjoint over~$k^a(\cB)$. Let~$\cC'$ (resp.~$\cC''$) be the smooth, projective, integral
curve corresponding to the compositum  of~$k^a(\cC)$ with~$k^a(\cB')$
(resp.~$k^a(\cB'')$). We have the following diagram:

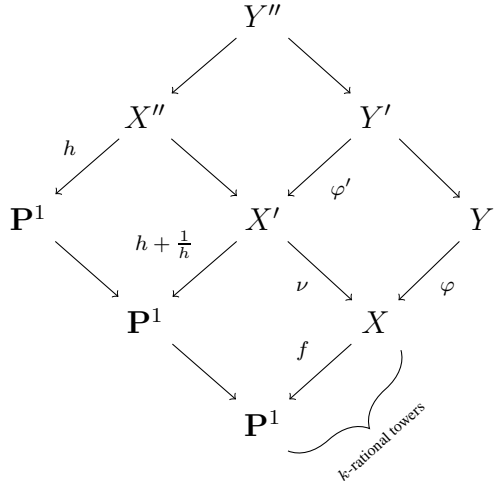
\begin{center}
\begin{tikzpicture}
\matrix[matrix of math nodes,row sep=0.75cm,column sep=0.75cm]
{
& & |(Ysec)|\cC'' & & \\
&|(Xsec)|\cB'' & & |(Yprim)|\cC' & \\
|(P1up)|\PP^1 & & |(Xprim)|\cB' & & |(Y)|\cC \\
& |(P1mid)|\PP^1 & & |(X)|\cB &  \\
& & |(P1down)|\PP^1 & & \\
};
\begin{scope}[every node/.style={pos=0.5,font=\scriptsize}]
\draw[->] (Xsec) -- (P1up) node[auto=right] {$h$} ;
\draw[->] (Xprim) -- (P1mid) node[auto=right] {$h+\frac{1}{h}$} ;
\draw[->] (X) -- (P1down) node[auto=right] {$f$} ;
\draw[->] (Ysec) -- (Yprim) ;
\draw[->] (Xsec) -- (Xprim) ;
\draw[->] (P1up) -- (P1mid) ;
\draw[->] (Ysec) -- (Xsec) ;
\draw[->] (Yprim) -- (Xprim) node[auto=left] {$\varphi'$} ;
\draw[->] (Y) -- (X) node[auto=left] {$\varphi$} ;
\draw[->] (Yprim) -- (Y) ;
\draw[->] (Xprim) -- (X) node[auto=right] {$\nu$} ;
\draw[->] (P1mid) -- (P1down) ;
\end{scope}
\draw[decorate,decoration={brace,raise=0.5cm,amplitude=5mm}]
(X.center) -- (P1down.center) node[below=1cm,pos=0.5,sloped] {\tiny $k$-rational towers};
\end{tikzpicture}
\end{center}
%\begin{equation*}
%\xymatrix@=20pt{
%&&\cC''\ar@{->}[dl]\ar@{->}[dr]_{\langle w\rangle}& & \\
%&\cB''  \ar@{->}[dl]_{h}  \ar@{->}[dr]_{\langle w\rangle} 
%&&\cC'  \ar@{->}[dl]_{\varphi'}\ar@{->}[dr] &\\
%\PP^1 \ar@{->}[dr]_{\langle w\rangle} &&\cB'  \ar@{->}[dl]_{h + \frac{1}{h}}  \ar@{->}[dr]^\nu
%&&\cC\ar@{->}[dl]_\varphi\\
%&\PP^1 \ar@{->}[dr] &&\ar@{->}[dl]_f\cB&  \\
%&&\PP^1 &}\label{diagramme:rev}
%\end{equation*}
The cover~$\cC'' \to \cC$ is again a $D_p$-Galois cover
and the cover~$\cC'' \to \cC'$ has degree~$2$.

Let us show that the cover~$\varphi' : \cC' \to \cB'$ has the
expected properties.

First of all, it is clear that the construction above yields a morphism of
stacks~$\FF : \Mod_\varphi \to \Mod_{\varphi'}$.
The Galois equivariance is a direct consequence of the fact that the middle diagonal tower
is defined over~$k$. This morphism is in fact
fully faithful because the sub-cover~$\varphi : \cC \rightarrow \cB$ 
of~$\nu\circ \varphi' : \cC'\rightarrow  \cB$  is
the maximal sub-cover unramified at~$f^{-1}(\{2,-2,\infty\})$.

Last we have to prove that the curve~$\cB'$ does not have any
non-trivial automorphism. Let~$\theta'$ be an automorphism of~$\cB'$. Call~$Z$ the image~$Z$
of~$\nu\times (\nu\circ \theta') : \cB'
\rightarrow (\cB\times\cB)$. Let~$\pi_1 : \cB \times \cB \to \cB$ be the projection to the
first factor. The map~$\nu$ factors as:
$$
\nu : X' \longrightarrow Z \overset{\pi_1}{\longrightarrow} X
$$
and it has prime degree~$p$. So~$Z$ is either isomorphic to~$X$ or birationaly equivalent
to~$X'$. In the latter case, the geometric genus of~$Z$ would
be $>\frac{1}{4}\deg (f)p\ge 1+2p(g(\cB)-1)+p^2$ by Hurwitz genus formula.
But the bi-degree  of~$Z$ is~$\le (p,p)$; so, by lemma~\ref{lem:adjonction}, its virtual arithmetic
genus is less than~$1+2p(g(\cB)-1)+p^2$. A contradiction.
Therefore~$Z$ is a correspondence of bi-degree~$(1,1)$ which defines an automorphism~$\theta$
of~$\cB$ such that~$\theta\circ \nu =\nu \circ \theta'$. Such an automorphism preserves
the ramification data of~$\nu$, the one of its Galois closure~$\cB'' \to \cB$ and also 
the one of the unique subcover
of degree~$2$ of the cover~$\cB'' \to \cB$. Since this last cover is exactly ramified
above~$f^{-1}(\{-2,2\})$,
we deduce that~$\theta =\Id$ and then that~$\theta'$ is a $k^a$-automorphism of the cover~$\nu$.
Since~$\nu$ does not have any non-trivial automorphism, necessarily~$\theta' = {\Id}$.
\end{preuve}

\section{Quasi-projective surfaces}\label{section:surfacesaffines}

Let~$k$ be a field of characteristic zero. In this section, we give a general process
which associates to each $k^a$-cover of curves, a smooth quasi-projective
integral $k^a$-surface with the same field  of moduli and fields of definition.

\begin{theoreme}\label{theoreme:quasi}
Let $X_k$ be a smooth, projective, geometrically integral curve over $k$ and
set $X=X_k\times_{\Spec(k)}\Spec(k^a)$.
Let $Y$ be a smooth, projective, integral curve over $k^a$ and
let $\varphi : Y \rightarrow X$ be a non-constant  morphism of $k^a$
curves. Then there exists a smooth quasi-projective integral $k^a$-surface having the same field
of moduli and the same fields of definition as~$\varphi$.
\end{theoreme}

First of all, by propositions~\ref{prop:grandgenre} and~\ref{prop:sansauto}, one can 
assume  that
the base curve~$\cB$ has genus greater than~$2$ and has no non-trivial
$k^a$-automorphism.

We consider the product~$\cB \times \cC$ of the two curves
and we denote by~$G(\varphi)$ the graph of $\varphi$ inside this product.
Let~$U$ be the open complementary set
 of~$G(\varphi)$ in~$\cB \times \cC$.

The surface we are looking for  is nothing else than the open set~$U$.
We  call it 
the {\em mark} of the cover~$\varphi : \cC \to \cB$ and we now prove that 
is has the same  field of moduli and the same 
fields of definition as $\varphi$.

\medbreak

We need two lemmas.

\begin{lem}\label{lem:morphismeempreinte}
Let~$l/k$ be a finite extension of~$k$ inside~$k^a$.
Let $X_k$ be a smooth, projective, geometrically integral $k$-curve.
Set $X=X_k\times_{\Spec(k)}\Spec(k^a)$ and assume that the genus of~$X$ is greater
than~$2$ and that $X$ has no non-trivial $k^a$-automorphism.
Let~$U_l$ and~$V_l$ be the marks of two non-trivial
geometrically integral $l$-covers~$\varphi_l : \cC_l \to \cB_l$ and~$\psi_l : \cD_l \to
\cB_l$, where~$\cB_l = X_k\times_{\Spec(k)}\Spec(l)$.

Then every morphism of covers between~$\varphi_l : \cC_l \to \cB_l$ and~$\psi_l : \cD_l \to \cB_l$
induces a morphism between the corresponding marks~$U_l$ and~$V_l$. Conversely,
every surjective $l$-morphism from~$U_l$ to~$V_l$  is equal to
$\Id \times \gamma_l$ where $\gamma_l :\cC_l \rightarrow \cD_l$ is a  $l$-morphism
between the covers~$\varphi_l : \cC_l \to \cB_l$ and~$\psi_l : \cD_l \to \cB_l$.
\end{lem}

\begin{preuve}
Recall that a $l$-morphism between the covers~$\cC_l \overset{\varphi_l}{\longrightarrow} \cB_l$
and~$\cD_l \overset{\psi_l}{\longrightarrow} \cB_l$
is a $l$-morphism
of $l$-curves~$\gamma_l :\cC_l\rightarrow\cD_l$  such that~$\psi_l\circ\gamma_l= \varphi_l$.
The product~$\Id\times \gamma_l : \cB_l\times\cC_l \to\cB_l\times\cD_l$ 
maps the graph of~$\varphi_l$ to the graph of~$\psi_l$ and also the
mark~$U_l$ to the mark~$V_l$.

Conversely, let~$\upsilon_l$ be a surjective $l$-morphism form~$U_l$ to~$V_l$.
We denote by~$\upsilon : U \to V$, $\varphi : \cC \to \cB$, $\psi : \cD \to \cB$ the base change
to~$k^a$ of~$\upsilon_l$, $\varphi_l$, $\psi_l$ respectively.

% Let~$\pi_1,\pi_2$
%denote the projections onto the components~$\cB,\cC,\cD$.

Let~$y$ be a
closed $k^a$-point of~$\cC$. Let~$\pi_2 : \cB \times \cD \to \cD$ be the projection on the second
factor. The restriction of~$\pi_2 \circ \upsilon$
to~$(\cB\times\{y\})\cap U$ is a constant function because the genus of~$\cB$ is
less than the one of~$\cD$. We denote by~$\gamma(y)$ this constant; this defines a
morphism~$\gamma : \cC \rightarrow \cD$ which cannot be constant since~$\upsilon$ is
 surjective. Let~$\pi_1 : \cB \times \cD \to \cB$ be the projection on the first
factor.
The  restriction of~$\pi_1\circ \upsilon$ to~$(\cB\times\{y\})\cap U$ is a
morphism~$\beta_y$ with values in~$\cB$. Let~$F\subset \cC$  the set of closed $k^a$-points
of~$\cC$ such that the morphism~$\beta_y$ is constant. This
 is a closed set; and a finite one 
because~$\upsilon$ is surjective. For a closed $k^a$-point~$y\not \in F$ the
morphism~$\beta_y$ induces an automorphism of~$\cB$, which is trivial since~$\cB$ does not
have any non-trivial automorphism. Thus we have~$\upsilon (x,y)=(x,\gamma(y))$ for every
closed $k^a$-point~$x$ on~$\cB$ and~$y$ on~$\cC$ with~$y \not \in F$ and~$(x,y)\in U$.
Let~$x$ be a closed $k^a$-point of~$\cB$. The restriction of~$\pi_1 \circ \upsilon$
to~$(\{x\}\times \cC)\cap U$ is constant and equal to~$x$ on the non-empty open
set~$(\{x\}\times (\cC - F))\cap U  $.  So it is a constant function. So~$F$ is empty 
and~$\upsilon$ is the restriction of~$\Id \times \gamma$ to~$U$. 
Thus~$\Id\times \gamma$ maps $U$ to~$V$ and therefore~$\psi\circ\gamma=\varphi$.
Moreover~$\gamma$ must be defined over~$l$ since~$\upsilon, U, V$ are defined over~$l$.
\end{preuve}

\begin{lem}\label{lem:produitdefi}
Let $X_k$ be a smooth, projective, geometrically integral $k$-curve.
Set $X=X_k\times_{\Spec(k)}\Spec(k^a)$ and assume that the genus of~$X$ is greater
than~$2$ and that $X$ has no non-trivial $k^a$-automorphism.
Let~$U$ be the mark of a non-constant $k^a$-cover~$\varphi : \cC \to \cB$, where~$\cC$
is a smooth, projective, integral $k^a$-curve. Then:
\begin{enumerate}
%\item the group of $k^a$-automorphisms of~$U$  is equal to the group of $k^a$-automorphisms
%of the cover~$\varphi : \cC \to \cB$;
\item $k$ is the field of moduli of~$U$ (in the category of quasi-projective varieties)
if and only if it is the field of moduli of the cover~$\varphi : \cC \to \cB$;
\item an algebraic extension of~$k$ is a field of definition of~$U$ 
if and only if it is the field of definition of the cover~$\varphi : \cC\to\cB$.
\end{enumerate}
\end{lem}

\begin{preuve}
It is clear that the construction of the mark from the cover commutes with base change.
This yields a morphism of stacks~$\FF : \Mod_\varphi \to \Mod_U$ which is
fully faithful according to lemma~\ref{lem:morphismeempreinte}.
The result follows by proposition~\ref{obstruction_transportor}. In particular, $\FF$
has an inverse functor~$\GG : \Mod_U \to \Mod_\varphi$.
\end{preuve}

\section{Proper normal surfaces}\label{section:surfacesnormalesprojectives}

In this section $k$ is a field of characteristic zero.
We start from  a cover of curves, having $k$ as field of moduli,
and we construct a proper normal integral
surface over $k^a$,  having the same field of moduli and the same
fields of definition as the original cover.

\begin{theoreme}\label{theoreme:projective}
Let~$\cB_k$ be a smooth, projective, geometrically integral curve over~$k$
and set~$\cB = \cB_k \times_{\Spec(k)} \Spec(k^a)$.
Let~$\cC$ be a smooth projective, integral
curve over $k^a$ and let~$\varphi : Y \to X$ be a cover.
Assume that~$k$ is the field of moduli of~$\varphi$.
Then, there exists a proper, normal and integral surface $\cS$
over $k^a$,
having  $k$ as field of moduli, and having the same fields of definition
as $\varphi$.
\end{theoreme}

The proof of this theorem is given in the rest of this section. The surface is constructed
as the cover of a surface~$\cB \times \cD$, strongly ramified along the graph of~$\psi$,
where~$\psi : \cD \to \cB$ is a suitably chosen cover of curves deduced from~$\varphi$.

\subsection{Construction of the  surface~$\cS$}
\label{section:constructionsurfacenormale}

The construction of the surface is divided in several steps.

\bigbreak

{\noindent\bfseries Step 0. Starting point.}

\medbreak

We keep notation and assumptions of theorem~\ref{theoreme:projective}. We denote
by~$g(X)$ the genus of~$X$ and by~$d$ the degree of the cover~$\varphi$.
 According to proposition~\ref{prop:sansauto}, we may assume that $g(X)$ is at least~$2$ and
that $\cB$ has no non-trivial automorphism over $k^a$.

\bigbreak

{\noindent\bfseries Step 1. A system of generators~$f_1, \ldots, f_I$ of the function
field~$k(\cB_k)$.}

\medbreak

We need to exhibit some $k$-{\bfseries rational} functions on~$\cB$.

\begin{lem}[The functions~$f_i$ on~$\cB$ and the primes~$p_i$] \label{fcts_f_i}
There exist~$I \in \NN^*$, some prime integers $p_1, \ldots, p_I > d$, and
functions~$f_1, \ldots, f_I \in k(\cB_k)$ satisfying the following conditions:
\begin{enumerate}
\item the functions~$(f_i)_{1 \leq i \leq I}$ generate the field~$k(\cB_k)$ over~$k$;
\item for every~$1 \leq i \leq I$ and every~$\lambda \in k^a$, none of the
functions~$f_i - \lambda$ is a $p_i$-th power in~$k^a(\cB)$;
\item let $\Pi = \prod_{i=1}^I p_i$ and let $M$ (resp.  $m$) be the maximum (resp. minimum)
among the degrees of the  $f_i$, then:
$$
\forall 1 \leq i \leq I, \qquad 1+2(g(\cB)-1)\Pi+\Pi^2 < m \leq \deg(f_i) \leq M.
$$
\end{enumerate}
\end{lem}

\begin{preuve}
We first choose a finite generating system   
$(h_j)_{1\le j\le J}$ 
of  $k(\cB_k)$ over $k$. We assume that
none of  the   $h_j$ is a power in $k^a(\cB)$.
We set~$I = 2J$ and let~$\Pi = \prod_{i = 1}^{I} p_i$  be the 
product of the first  $I$ prime integers greater 
than the degree  $d$ of
$\varphi$. We choose two distinct prime integers $a$ and  $b$, both 
bigger than
$1+2(g(\cB)-1)\Pi+\Pi^2$.
For every  $1\le j\le J$, we set:
$$
f_j=h_j^a,
\qquad \text{ and } \qquad
f_{j+J}=h_j^b.
$$
We can choose~$a$ and~$b$ in such  a way that 
none of the functions~$f_i - \lambda$ is
a $p_i$-th power in $k^a(\cB)$ for~$\lambda \in k^a$ and~$1 \leq i \leq I$:
this is evident for  $\lambda=0$. If  $\lambda\not =0$
and if  $h_i^a-\lambda=\prod_{0\le k \le a-1}
(h_i-\zeta_{a}^k\lambda^{\frac{1}{a}})$ is a power,  then 
$h_i$ has at least   $a$ distinct singular values. This is impossible
if we choose an  $a$ bigger than the number
of singular values of  $h_i$.

We note also that the  $(f_i)_{1\le i\le I}$ generate
$k(\cB_k)$ over $k$ and that they all have a degree greater than
$1+2(g(\cB)-1)\Pi+\Pi^2$, as expected.
\end{preuve}

\bigbreak

{\noindent\bfseries Step~2. A cover~$\psi : \cD \to \cB$ of large enough degree.}

\medbreak

Let  $p$ be a prime integer bigger than
$(g(\cB)+IM)\Pi$. We call  $\cD$ the  curve
and  $\psi : \cD \rightarrow \cB$ the degree $pd$ cover
given by  proposition~\ref{prop:granddegre}.
The genus of  $\cD$ is bigger than  $dp>(g(\cB)+IM)\Pi$ and the
covers~$\varphi$ and~$\psi$ have the same field of moduli and the
same fields of definition.

\bigbreak

{\noindent\bfseries Step~3. A system of functions~$g_1, \ldots, g_I$ on~$\cB \times \cD$.}

\medbreak

Using the previous functions~$f_i$, we define functions on~$\cB \times \cD$.

\begin{lem}[The functions~$g_i$ on~$\cB \times \cD$] \label{fcts_g_i}
For every~$1 \leq i \leq I$,
let~$g_i$ be the function on $\cB\times\cD$ defined by:
$$
g_i (P,Q)=f_i(\psi(Q))-f_i(P).
$$
Then:
\begin{enumerate}
\item the negative part  $(g_i)_\infty$ of the divisor of  $g_i$
is~$(f_i)_\infty\times \cD+\cB\times (f_i\circ\psi)_\infty$;
\item the  positive parts   $(g_i)_0$ are such that~$\pgcd_i((g_i)_0)=G(\psi)$,
where  $G(\psi)$ is the  graph of $\psi$;
\item for every point~$P \in \cB$ the function~$Q \mapsto g_i(P,Q)$ on~$P \times \cD$ is not
a $p_i$-th power.
\end{enumerate}
\end{lem}

\begin{preuve}
The first two points are easy. To prove the third one, recall that each function~$f_i$
is such that none of the~$f_i-\lambda$ for~$\lambda \in k^a$ is a $p_i$-th power
(lemma~\ref{fcts_f_i}). Since the
degree~$pd$ of~$\psi$ is prime to~$p_i$, none of the function~$f_i\circ \psi - \lambda$ is a 
$p_i$-th power in $k^a(\cD)$. Condition~3 follows.
\end{preuve}

Let us note that, if~$\psi$ is defined over a field~$l$, then so are the functions~$g_i$.

\bigbreak

{\noindent\bfseries Step~4. At last, the surface~$\cS$.}

\medbreak

Let $k^a(\cB\times\cD)$ be
the field of functions of the surface $\cB\times\cD$.
We define a regular radicial  extension of
$k^a(\cB\times\cD)$ by setting
$$
y_i^{p_i}=g_i.
$$
We denote by  $\cS$ the  normalization of
 $\cB\times \cD$ in the latter radicial extension.
It is a normal surface by construction and there is a ramified cover:
$$
\chi : \cS\rightarrow \cB\times \cD
$$
which is a Galois cover of surfaces over~$k^a$ with Galois group $\prod_{i=1}^{I}\ZZ/p_i\ZZ$.

\subsection{The group of automorphisms of~$\cS$}

We denote by~$A$ the  group of  $k^a$-automorphisms of~$\psi$. 
An element in $A$ induces   
a  $k^a$-automorphism of the surface~$\cB\times\cD$, and this latter automorphism  can be
lifted uniquely to  
an  automorphism of $k^a(\cS)/k^a$ that
fixes all  $y_i$ and stabilizes  $k^a( \cB\times \cD)$.
In the sequel we shall use the same notation for an automorphism
of~$\psi$, the induced automorphism of~$\cB\times\cD$ and
its lift to~$\cS$.
In other words, $A$ can be identified with  a subgroup
of~$\Aut_{k^a}(\cS)$, the  group of $k^a$-automorphisms of~$\cS$.

We know another subgroup of~$\Aut_{k^a}(\cS)$, namely the Galois 
group~$B=\prod_{i=1}^I \ZZ/p_i\ZZ$ of the extension $k^a(\cS)/k^a(\cB \times\cD)$.

To summarize, $A$ is the set of~$\alpha$ such that the following diagram commute:
\begin{equation}\label{diag:Aplusieurspointsdevue}
\vcenter{
\xymatrix@=15pt{
\cS \ar@{->}[r]^{\alpha} \ar@{->}[d]_{\chi} & \cS \ar@{->}[d]^{\chi} \\
\cB \times \cD \ar@{->}[r]^{\alpha} & 
\cB \times \cD
}
}
\end{equation}
while $B$ is the set of~$\beta$ such that the following diagram commute:
\begin{equation}\label{diag:Bplusieurspointsdevue}
\vcenter{
\xymatrix@=15pt{
\cS \ar@{->}[rr]^{\beta} \ar@{->}[rd]_{\chi} & & \cS \ar@{->}[ld]^{\chi} \\
& \cB \times \cD
}
}
\end{equation}

It is clear that~$A \times B \subset \Aut_{k^a}(\cS)$. We now prove
that this inclusion is an equality. To this end, we introduce a family of curves
on $\cS$.

\begin{lem}[The curves~$\cE_Q$] \label{curves_E_Q}
For any point~$Q$  on $\cD$,
we call~$\cE_Q$   the
inverse image of $\cB \times Q$ by  $\chi$
and we denote by~$\chi_Q : \cE_Q\rightarrow \cB \times Q$ the restriction
of $\chi$ to $\cE_Q$.
The geometric genus of
 $\cE_Q$ can be bounded from above:
\begin{equation} \label{majoration_genre} %\tag{$\star$}
g(\cE_Q) \leq( g(\cB)+IM)\Pi < g(\cD),
\end{equation}
and the genus of any non-trivial subcover of~$\chi_Q$ can be bounded from below:
\begin{equation} \label{minoration_genre} %\tag{$\star\star$}
1+2(g(\cB)-1)\Pi+\Pi^2 < m \leq g(\text{non-trivial subcover of~$\chi_Q : \cE_Q \to \cB$}).
\end{equation}
\end{lem}

\begin{preuve}
If $Q$ is the generic   point
on  $\cD$, then  $\cE_Q$ is a geometrically integral
curve and~$\chi_Q$ is a degree $\Pi$, geometrically connected cover.
The degree of the ramification divisor
of this cover is bounded from above by the product
$2IM$ (where  $I$ is the number 
of functions in the family  $(f_i)_i$
and  $M$ is the maximum of the degrees of these functions).
The upper bound follows.

For the lower bound, let us consider a non-trivial subcover of ~$\chi_Q$. Such a cover has degree
at least~$p_1 \geq 3$ and its ramification divisor has degree at least~$m$
(where~$m$ is the minimum among the degrees of the functions~$f_i$). So its genus is greater
than~$m$ and the lower bound follows.
\end{preuve}

\begin{lem}\label{lem:automorphismes}
The  group~$\Aut(\cS)$ of  $k^a$-automorphisms
of  $\cS$ is $A\times B$.
\end{lem}

\begin{preuve}
Let~$\theta$ be a  $k^a$-automorphism  of $\cS$.

First of all, let~$Q$ be the generic point of~$\cD$.
We know from 
inequality~(\ref{majoration_genre}) of lemma~\ref{curves_E_Q} that $g(\cE_Q) < g(\cD)$.
We deduce that~$\theta(\cE_Q)=\cE_{\alpha(Q)}$
where  $\alpha$ is a $k^a$-automorphism of $\cD$.

%Montrons que~$\beta \overset{\text{d{\'e}f.}}{=} \theta \circ \alpha^{-1} \in B$, c'est-{\`a}-direque~$\chi \circ \theta \circ \alpha^{-1} = \chi$ ouencore~$\chi \circ \theta = \chi \circ \alpha = \alpha \circ \chi$. Il suffit de montrer cette{\'e}galit{\'e} sur la fibre g{\'e}n{\'e}rique~$\cE_Q$ pr{\'e}c{\'e}dente. Notons~$\chi_Q : \cE_Q \to \cB$et~$\chi_{\alpha(Q)} : \cE_{\alpha(Q)} \to \cB$ les projections respectives sur~$\cB$ via~$\chi$.

We now prove that
the isomorphism between~$\cE_Q$ and $\cE_{\alpha(Q)}$ induced  by~$\theta$
makes the following diagram  commute:
\begin{equation} \label{diag:E_Q_E_alphaQ}
\vcenter{
\xymatrix@=30pt{
\cE_Q \ar@{->}[r]^{\theta} \ar@{->}[d]_{\chi_Q} & \cE_{\alpha(Q)} \ar@{->}[d]^{\chi_{\alpha(Q)}} \\
\cB \times Q \ar@{->}[r]^{ \Id \times \alpha} & \cB \times \alpha(Q)
}
}
\end{equation}
%avec~${\Id} \times \alpha$ comme isomorphisme en bas.
Indeed, the cartesian product of the maps  $\chi_Q$ and  
$\chi_{\alpha (Q)}\circ \theta$ defines a morphism:
$$
\cE_Q 
\xrightarrow{\chi_Q \times \left(\chi_{\alpha (Q)}\circ \theta \right)}\cB \times \cB,
$$
whose image $W$
is a divisor with bidegree $\le (\Pi,\Pi)$. Using
lemma \ref{lem:adjonction} we deduce that the arithmetic genus of
$W$  is smaller than or equal to
 $1+2(g(\cB)-1)\Pi+\Pi^2$. Let~$\pi_1 : X \times X \to X$ be the projection on the first
factor. The morphism~$\chi_Q$ factors as:
$$
\chi_Q :  \cE_Q \longrightarrow W \overset{\pi_1}{\longrightarrow} X.
$$
The map $W \overset{\pi_1}{\longrightarrow} X$ is a birational isomorphism. Otherwise, it would
define a non-trivial sub-cover of $\chi_Q :\cE_Q\rightarrow \cB$.
But we know 
from inequality~(\ref{minoration_genre}) of
lemma~\ref{curves_E_Q} that 
such a subcover  has geometric genus greater than or equal to 
$m>1+2(g(\cB)-1)\Pi+\Pi^2$. A contradiction. We deduce
that $W$ is a correspondence
of bidegree  $(1,1)$. Since $\cB$ has no non-trivial
 $k^a$-automorphism we deduce that
diagram~$(\ref{diag:E_Q_E_alphaQ})$ commutes.

We now prove that~$\alpha \in A$. We just showed that~$\theta$ induces an 
isomorphism between the covers~$\chi_Q : \cE_Q \to \cB$ 
and~$\chi_{\alpha(Q)} : \cE_{\alpha(Q)} \to \cB$.
Therefore these two covers have the same ramification
data: for every $1\le i \le I$ the points~$P$ such that~$f_i(P) = f_i(\psi(Q))$ and those such that~$f_i(P) = f_i(\psi(\alpha(Q)))$ are the same.
Thus:
$$
\forall i, \quad f_i(\psi(Q)) = f_i(\psi(\alpha(Q)))
$$
\noindent therefore $\psi(Q) = \psi(\alpha(Q))$, 
because the~$f_i$ 
generate  $k^a(\cB)$ over $k^a$ (lemma~\ref{fcts_f_i}).  So~$\psi = \psi \circ \alpha$,
and~$\alpha\in A$.

Diagram~\ref{diag:E_Q_E_alphaQ} implies
that the map~$\chi_{\alpha(Q)}\circ \theta : \cE_Q \to \cE_{\alpha(Q)}$
is equal to~$(\Id \times \alpha ) \circ \chi_Q$. And
this is~$\chi_{\alpha(Q)}\circ \alpha$ according to
diagram~$(\ref{diag:Aplusieurspointsdevue})$.
We set    $\beta= \theta\circ \alpha^{-1}$ and we check that
$\chi_{\alpha(Q)}\circ
\beta=\chi_{\alpha(Q)}$. Since 
 $Q$ is generic and  $\alpha$ surjective
we deduce  that $\chi\circ
\beta=\chi$ so  $\beta \in  B$.
We conclude that~$\theta=\beta \alpha \in A\times B$ as was to be shown.
\end{preuve}

\begin{remarque}
We have proven something slightly stronger than
lemma \ref{lem:automorphismes}: the group of birational
$k^a$-automorphisms
of  $\cS$ is  $A\times B$. We shall not need this stronger result.
\end{remarque}

\subsection{Field of moduli and fields of definition of~$\cS$} \label{s_def_S}

To prove theorem~\ref{theoreme:projective}, we have to show that the cover~$\varphi$
and the surface~$S$ share the same field of moduli and the same
fields of definition. In fact, by construction, one can replace the cover~$\varphi$
by the cover~$\psi$, since those two covers have same field of moduli and fields of definition.

The construction of section~\ref{section:constructionsurfacenormale} yields a
morphism of stacks~$\FF : \Mod_\psi \to \Mod_S$. To see this, let us consider~$l/k$
extension inside~$k^a$ and let~$\psi_l : \cD_l \to \cB_l$ be an $l$-model of~$\psi$.
We just follow the line of the construction, replacing~$\psi$ by~$\psi_l$.
Since the functions~$f_i$ are $k$-rational, the functions~$g_i$ lie in~$l(\cB_l\times\cD_l)$.
Then the radical extension defined by the equations~$y_i^{p_i} = g_i$ is a regular
extension of~$l(\cB_l\times\cD_l)$. The normalization of~$\cB_l\times\cD_l$ in this extension
is a surface~$\cS_l$ which is defined over~$l$. Of course, this surface~$\cS_l$ is an $l$-model
of~$\cS$ and the morphism~$\FF$ is defined on objects by~$\FF(\psi_l) = \cS_l$. Because
functions~$f_i$ are $k$-rational,~$\FF$ is a morphism of stacks.
By proposition~\ref{obstruction_transportor},
$k$ is the field of moduli of~$\cS$ and every field of definition of~$\psi$ (or~$\varphi$)
is a field of definition of~$\cS$.

Unfortunately,~$\FF$ is not fully faithful. As in
proposition~\ref{prop:granddegre}, we use proposition~\ref{prop:quotient_by_subgroup}
to construct a morphism the other way around. The group~$\Aut_{k^a}(\cS)$ is a normal subgroup
 of~$\Aut_k(\cS)$. Conjugation by an element of~$\Aut_k(\cS)$ induces an automorphism
of~$\Aut_{k^a}(\cS)$. This automorphism must stabilize the unique sub-group of order~$\Pi$
of~$\Aut_{k^a}(\cS)$, which is nothing but~$B = \Aut_{k^a}(\chi)$.
Let~$U$ be the mark of the cover~$\psi$. This is the complementary set of the branch locus
of the quotient map~$\chi : S \to X \times Z$. According to
proposition~\ref{prop:quotient_by_subgroup},
taking the complementary set of the branch locus of a quotient map defines
a morphism of
stacks $\GG : \Mod_S \to \Mod_U$.
Therefore, every field of definition of~$S$ is a field of definition of the mark of~$\psi$
and then also a field of definition of~$\psi$ by lemma~\ref{lem:produitdefi}. Indeed
the proof of this lemma provides a morphism from~$\Mod_U$ to~$\Mod_\psi$ and the proof
of proposition~\ref{prop:granddegre} provides a morphism from~$\Mod_\psi$ to~$\Mod_\varphi$.

\section{Curves}\label{section:courbes}

In this  section $k$ is still a field of characteristic zero.
We start from  a cover of curves, having $k$ as field of moduli,
and we construct, a projective normal integral
curve  over $k^a$,  having the same field of moduli and the same
fields of definition as the original cover. This will prove theorem~\ref{theoreme:courbe}.

%Le rev{\^e}tement de degr{\'e} $\cE\rightarrow \cC$ n'est pas $k^a$-galoisien donc $\varphi$ est  le quotient de $\psi$ par son groupe d'automorphismes

We shall make use of the surface~$\cS$ constructed
in section \ref{section:surfacesnormalesprojectives}. 
So we keep the  notation of
 section \ref{section:surfacesnormalesprojectives}. 
We know that  $\cS$
has field of moduli $k$ and the same fields of definition
as the initial cover
$\varphi : \cC \rightarrow \cB$
(or equivalently 
 $\psi : \cD \rightarrow \cB$).

The main idea is to draw on~$\cS$ a singular (but stable)
curve inheriting the field  of moduli and fields of
definition of~$\cS$; then to deform it
to obtain a smooth projective curve.

\subsection{Two  stable curves}\label{subsection:deux}

In section~\ref{section:constructionsurfacenormale},
we have constructed a cover
$\chi : \cS\rightarrow \cB\times\cD$
strongly ramified along the graph
of  $\psi : Z \to X$. 
 For any point~$P$  on $\cB$,
we call~$\cF_P$   the
inverse image of $P \times \cD$ by  $\chi$ and
  $\chi_P : 
\cF_P \rightarrow P\times \cD$   the   correstriction
of  $\chi$ to  $P\times \cD$.
We call~$\Gamma$
the union of the supports of all
divisors of the  functions~$g_i$ of lemma~\ref{fcts_g_i}. It contains the
ramification locus of the cover~$\chi$.

%JMC

\begin{lem} \label{fcts_fg}
There exist two  non-constant $k$-rational functions~$f,g \in k^a(\cB)$ 
such that:
\begin{enumerate}
\item\label{fcts_fg_Zf} the divisor~$((f)_0 + (f)_\infty) \times \cD$ crosses
transversally~$\Gamma$;
\item\label{fcts_fg_Zg} the divisor~$\cB \times ((g \circ
  \psi)_0+(g\circ \psi)_\infty)$ crosses transversally~$\Gamma \cup [((f)_0 +(f)_\infty )\times \cD]$;
%JMC je coupe en deux
\item\label{fcts_fg_Aut_gpsi} any  $k^a$-automorphism of~$\cD$ that stabilizes
the fiber~$(g \circ \psi)_0$
is an  automorphism of the cover~$\psi$ (note that the preceding
condition implies that this fiber
is simple);
%JMC ci dessous j'enl{\`e}ve le mot ``point''
\item\label{fcts_fg_Aut_gpsichi} for any zero  $P$ of  $f$, 
the cover 
$g\circ \psi \circ \chi_P  : \cF_P \rightarrow \PP^1$ has no
 automorphism other  than the elements of  $A\times B$:
$$
\Aut_{k^a}(g\circ \psi \circ \chi_P)=\Aut_{k^a}(\psi\circ \chi_P)=A\times B.
$$
\end{enumerate}
\end{lem}

\begin{preuve}
Let  $f\in k^a(\cB)$ be a $k$-rational non-constant function.
We apply  lemma \ref{lem:fibrescoupenttransversalement} 
to  $k$, $\cB$,   $\cD$, $\Gamma$ and $f$. We deduce that there exist
two distinct scalars $x$ and  $y$ in $k^a$ such that  $(f)_x\times \cD$ and
 $(f)_y\times
\cD$ cross transversally  $\Gamma$. We even can choose~$x$ and~$y$ in $k$
and such that for every point~$P$ in~$f^{-1}(x)$ or~$f^{-1}(y)$, the fiber
of every function~$f_i \in k(\cB)$ above~$f_i(P)$ does no meet
the singular values of~$\psi$, that is:
\begin{equation} \label{contraintes_valeurs_singulieres}
\forall P \in f^{-1}(x) \cup f^{-1}(y),
\qquad f_i^{-1}(f_i(P)) \cap \{\text{singular values of~$\psi$}\} = \emptyset.
\end{equation}
We replace  $f$ by  $(f-x)/(f-y)$ and the first condition
is fulfilled.

Now, for every zero  $P$  of $f$,  we see  that  $\cF_P$ is 
smooth and geometrically integral, because $(f)_0\times \cD$ crosses
transversally the ramification locus  $\Gamma$ of $\chi$.  
We now prove that 
$$
\Aut_{k^a}(\psi\circ\chi_P)=A\times B.
$$
Indeed the function field~$k^a(\cF_P)$
is the  compositum:
$$
\xymatrix@=15pt{
& k^a(\cF_P) \ar@{-}[dl] \ar@{-}[dr] \ar@{-}[dd] \\
k^a(\cD) \ar@{-}[dr]_{pd} & & k^a(\cB') \overset{\text{def.}}{=}
k^a\left( (f_i - f_i(P))^{\frac{1}{p_i}}, 1\leq i \leq I\right) \ar@{-}[dl]^{\prod_i p_i} \\
& k^a(\cB)
}
$$
where~$\cB' \to \cB$ is an abelian
cover with Galois group~$B = \prod_{i=1}^I \ZZ/p_i\ZZ$. 
The $k^a(\cB)$-extensions $k^a(\cD)$
and~$k^a(\cB')$ are linearly  disjoint (their degrees are coprime
and one of them is Galois) and condition~(\ref{contraintes_valeurs_singulieres}) implies that the extension $k^a(\cD)/k^a(\cB)$  is not ramified
above the zeros of the  functions $f_i-f_i(P)$.

Now,  any subcover of~$\cB' \to \cB$ is ramified above
the  zeros of at least one of the  functions~$f_i - f_i(P)$.
The same is true for any subcover of~$\cF_P \to \cD$. 
We deduce that~$\cD \to \cB$ is the maximal subcover
of~$\cF_P \to \cB$ that is not ramified above the
zeros of the functions~$f_i - f_i(P)$. Therefore any 
 $k^a(\cB)$-automorphism of~$k^a(\cF_P)$
stabilizes~$k^a(\cD)$. Thus:
$$
\Aut_{k^a(\cB)}(k^a(\cF_P)) = \Aut_{k^a(\cB)}(k^a(\cD)) \times \Aut_{k^a(\cB)}(k^a(\cB')),
$$
as was to be shown.

Next we look for a function  $g$ in $k(\cB)$ such that  $g\circ \psi$
has no  $k^a$-automorphism  but elements of  $A$ and, 
for every zero  $P$ of $f$, the cover $g\circ \psi \circ \chi_P$ has no  $k^a$-automorphism
but elements of  $\Aut_{k^a}(\psi\circ\chi_P)=A\times B$. 
According to  lemma
\ref{lem:fonctionsauto}, the functions in $k(\cB)$ that do not
fulfill all these conditions lie in a finite union of strict
sub-$k$-algebras.
Therefore there exists such  a function $g$. 

According to  lemma
\ref{lem:fibrescoupenttransversalement},
the scalars $x$ in $k$ such that  $(g\circ\psi)_x$ does not cross
$\Gamma \cup [((f)_0 +(f)_\infty ) \times \cD]$ transversally, are
finitely many.

According to lemma 
\ref{lem:autofibres},
the  $x$ in  $k$ such that $(g\circ\psi)_x$ has
a  $k^a$-automorphism not in  $A$, are 
finitely many.

Therefore there exist   two distinct scalars
$x$ and  $y$ in  $k$ such that  
$(g\circ\psi)_x$  and  $(g\circ\psi)_y$  cross
$\Gamma \cup [((f)_0 + (f)_\infty )\times \cD]$ transversally and 
$(g\circ\psi)_x$ has no  automorphism but those
in   $A$.
We replace  $g$ by  $(g-x)/(g-y)$ and the last three
 conditions are satisfied.
\end{preuve}

\bigbreak

{\bfseries\noindent The curves~$\cJ_0$ and~$\cK_0$.}

Let~$\cJ_0$ be the curve on~$\cB \times \cD$ with equation:
$$
f(P) \times g\circ \psi(Q) = 0.
$$
Let~$\cK_0$ be the inverse image of~$\cJ_0$ by~$\chi$.
These are singular curves over $k^a$. The two following lemmas are concerned with
the stability and the automorphism groups of these two curves.

\begin{lem}\label{lem:automorphismes2}
The curve~$\cJ_0$ is  stable and~$\Aut_{k^a}(\cJ_0) \simeq A$.
\end{lem}

\begin{preuve}
The curve~$\cJ_0$ is geometrically reduced because
the zeros of~$f$ and~$g\circ \psi$ are simple.
The singular points on $\cJ_0$ are the
couples~$(P,Q)$ on 
$\cB \times \cD$ such that~$f(P) = g\circ \psi(Q) = 0$. These are ordinary double  points.
Therefore $\cJ_0$  is semistable.
It is geometrically connected also.
Its irreducible  components
are isomorphic to~$\cB$ or~$\cD$. So they all
have genus~$\geq 2$. Therefore~$\cJ_0$ is a  stable curve.

We now prove that~$\Aut_{k^a} (\cJ_0) \simeq A$, i.e. that the
group of  $k^a$-automorphisms of  $\cJ_0$ is the group  $A$ of
$k^a$-automorphisms of $\psi$. It is clear that $A$ is included
in~$\Aut_{k^a} (\cJ_0)$. Conversely, 
let  $\theta$ be a $k^a$-automorphism of
 $\cJ_0$. Then $\theta$ permutes the irreducible 
components of $\cJ_0$. Some of these components are
isomorphic to  $\cB$, and the others are isomorphic to  
$\cD$.
Since  $\cB$ and $\cD$ are not $k^a$-isomorphic,
$\theta$ stabilizes the two subsets of components.
  If we restrict $\theta$ to
a component  isomorphic to  $\cB$ then  compose with
the projection on  $\cB$, we obtain a non-constant
$k^a$-morphism from $\cB$ to itself.
This morphism must be the identity because 
$\cB$ has no non-trivial $k^a$-automorphism. 
Therefore  $\theta$ stabilizes  each component
 isomorphic to  $\cD$. 
The singular points on such a component
are the zeros of
 $g\circ \psi$. 
%JMC ici je parle d'ensemble
The set of these zeroes
is  stabilized by no $k^a$-automorphism 
of $\cD$ but those of $\psi$ by~(\ref{fcts_fg_Aut_gpsi}) of lemma~\ref{fcts_fg}.
So the restriction of 
$\theta$ to any component isomorphic to $\cD$
is an element in $A$. If we compose $\theta$ with some well chosen
element in $A$, we may assume that
$\theta$ is trivial on one component isomorphic
to  $\cD$.
Therefore $\theta$ stabilizes every component 
isomorphic to  $\cB$. 
Since these components have no non-trivial automorphism, 
$\theta$ acts trivially on them.
Now let  $P\times \cD$ be a 
 component of  $\cJ_0$ isomorphic to  $\cD$. 
The restriction of $\theta$ to it
is an automorphism that fixes the singular points.
 These  points are the zeros of $g\circ \psi$. So
 $\theta$ restricted to $P\times \cD$
is an element of  $A$.  But  $A$ acts faithfully
on the set of zeros of $g\circ \psi$. We deduce
that  $\theta$  acts trivially 
on every component isomorphic to $\cD$.
\end{preuve}

Controlling the full group of $k^a$-automorphisms of~$D_0$
seems difficult to us.
So we shall be interested in the subgroup consisting of
{\it admissible} automorphisms. This subgroup
is denoted~$\Aut_{k^a}^{\text{adm.}} (\cK_0)$. 
We now explain what we mean by an admissible
automorphism.

We first notice that the
components of  $\cK_0$ are of two different kinds.
Some of them are covers of  some $\cB\times Q$ where
 $Q$ is a  zero of  $g\circ\psi$. We denote such 
a component by  $\cE_Q$. The other components 
are covers of some  $P\times \cD$ where
$P$ is a  $k^a$-zero of $f$. Such a component
is denoted by $\cF_P$. We call  $\chi_P : \cF_P\rightarrow P\times \cD$ and
$\chi_Q : \cE_Q \rightarrow \cB\times Q$ the
 restrictions of $\chi$ to components of
 $\cK_0$. Now let   $T$  be a singular point
on $\cK_0$ such that  $\chi(T)=(P,Q)$. So  $T$ lies in the intersection
of  $\cE_Q$ and  $\cF_P$. The  point on $\cE_Q$ corresponding to
 $T$ is denoted $U$. The point on  $\cF_P$ corresponding to  $T$
is denoted $V$. So  $\chi_Q(U)=P$ and  $\chi_P(V)=Q$.
And   $f\circ \chi_Q$ is a uniformizing
parameter for  $\cE_Q$ at  $U$, while
 $g\circ \psi \circ \chi_P$ is a uniformizing parameter for
 $\cF_P$ at $V$. Let   $\theta$ be  an  automorphism of  $\cK_0$ and
let  $T'=(U',V')$ be  the image of   $T=(U,V)$ by
$\theta$. We write   $\chi(T')=(P',Q')$.
We observe that  $f\circ \chi_{Q'}\circ \theta$ is a  uniformizing 
parameter for $\cE_Q$ at  $U$ and  $g\circ \psi \circ \chi_{P'}\circ \theta$
 is a uniformizing parameter for  $\cF_P$ at $V$.

We say that $\theta $ is an 
{\it  admissible automorphism}
of  $\cK_0$ if for every  singular   point $T$
 of 
$\cK_0$ we have 
\begin{equation}\label{eq:admissible}
\frac{f\circ \chi_{Q'}\circ \theta}{f\circ \chi_Q}(U)\times 
\frac{g\circ \psi\circ \chi_{P'} \circ \theta}{g\circ \psi  \circ \chi_P}(V)=1
\end{equation}
where  $\chi (T)=(P,Q)$ and  $\chi(\theta(T))=(P',Q')$.

The justification for this definition  is given
at  paragraph \ref{subsection:deformation}. Admissible
 automorphisms form a subgroup of the group
of $k^a$-automorphisms of $\cK_0$.

\begin{lem}\label{lem:sa}
The curve~$\cK_0$ is stable and~$\Aut_{k^a}^{\text{adm.}}(\cK_0) \simeq A \times B$.
\end{lem}

\begin{preuve}
First, it is clear that 
$A\times B$ acts faithfully on $\cK_0$,  and the corresponding  automorphisms 
are  admissible.

The curve~$\cK_0$ is drawn  on~$\cS$. Let us prove that
$\cK_0$ is a stable curve.
Points~(\ref{fcts_fg_Zf}) and~(\ref{fcts_fg_Zg}) of lemma~\ref{fcts_fg} imply that
the ramification locus~$\Gamma$ of~$\chi$ does not contain any singular points of~$\cJ_0$.
Therefore every 
singular point on~$\cJ_0$ gives rise to $\deg(\chi)$ singular points on~$\cK_0$;
and all these singular points are ordinary double points.
To prove that~$\cK_0$ is connected, we observe that the function~$g_i$ restricted
to any irreducible component of~$\cJ_0$
is not a $p_i$-th power because none of the functions~$f_i - \lambda$, $\lambda \in k^a$ is a
$p_i$-th power (and the  $f_i\circ \psi -\lambda$ are not either) as shown
in lemma~\ref{fcts_f_i}. Also the irreducible components
of $\cK_0$ correspond bijectively to those of $\cJ_0$.

Now let us prove that~$\Aut_{k^a}^{\text{adm.}}(\cK_0) \simeq A \times B$.
The components   $\cF_P$ and  $\cE_Q$ have different genera. Therefore no
$\cF_P$ is $k^a$-isomorphic to some  $\cE_Q$.
Thus any  $k^a$-automorphism $\theta $ of $\cK_0$ stabilizes  the set of all components
$\cF_P$ (and also the set of all $\cE_Q$).

Let~$Q$ and~$Q'$ be two  $k^a$-zeros of~$g \circ \psi$ such that~$\theta(\cE_Q) = \cE_{Q'}$. 
As in the proof of lemma~\ref{lem:automorphismes}, we notice that the image of~$\cE_Q$ in the product~$\cB \times \cB$, 
by the morphism~$\chi_Q \times \chi_{Q'} \circ \theta$, has an arithmetic genus smaller than or equal to ~$1+2(g(\cB)-1)\Pi+\Pi^2$. 
Again, this implies that this  image is
 $k^a$-isomorphic to~$\cB$ (otherwise, this image would have geometric genus bigger than~$1+2(g(\cB)-1)\Pi+\Pi^2$ by Hurwitz formula). 
Since~$\cB$ has no $k^a$-automorphism, we deduce as before that~$\theta$ induces an isomorphism of covers between the
 restrictions~$\chi_Q : \cE_Q \to \cB$ and~$\chi_{Q'} : \cE_{Q'} \to \cB$ of~$\chi$.
Thus
\begin{equation}
\chi_Q = \chi_{Q'} \circ \theta.\label{eq:chiq}
\end{equation}
This implies that~$\theta$ stabilizes every  component~$\cF_P$ where~$P$ is any $k^a$-zero of~$f$. Indeed,
let us start from a singular  point~$T =(U,V)\in  \cE_Q \cap \cF_P $,
where~$P$ is a $k^a$-zero of~$f$
and~$Q$ is a  $k^a$-zero of~$g \circ \psi$.
We thus have~$\chi(T) = (P,Q) \in \cB \times \cD$.
We know there  exists~$P' \in \cB(k^a)$ and~$Q' \in \cD(k^a)$ such
that~$\theta(T) \in \cF_{P'} \cap \cE_{Q'}$.
%So~$\theta(T) \in \cE_{Q'}\cap \theta(\cE_{Q})$ and~$\theta(\cE_{Q}) = \cE_{Q'}$.
We deduce from Equation~(\ref{eq:chiq}) that:
$$
P' = \chi_{Q'} \circ \theta(T) = \chi_Q(T) = P. 
$$
We conclude that~$P = P'$ and~$\theta(\cF_P) = \cF_P$. 

Now, we deduce from formulae~(\ref{eq:admissible}) and (\ref{eq:chiq}) that:
\begin{equation} \label{gpsi_Q'_sur_gpsi_Q}
\frac{g \circ \psi \circ \chi_P \circ \theta}{g \circ \psi \circ \chi_P}(V) = 1.
\end{equation}

Call~$\theta_P$ the restriction of~$\theta$ to~$\cF_P$. This is an  automorphism of~$\cF_P$. We prove that 
$\theta_P$ is the  restriction to~$\cF_P$ of an  element of~$A \times B$. To this end, we introduce the function~$h_P = g \circ \psi \circ \chi_P \in k^a(\cF_P)$. 
The degree of $h_P$ is~$\deg(g) \times pd \times \Pi$ and its  zeros are all simple. These zeros are the
 intersection points between~$\cF_P$ and the other components of~$\cK_0$. Since~$\theta_P$ permutes
these zeros, the functions~$h_P \circ\, \theta_P$ and~$h_P$ have the same divisor of zeros.
Therefore the only possible  poles of the function~$\frac{h_P}{h_P \circ\, \theta_P} - 1$ are the poles of~$h_P$. 
Thus the degree of ~$\frac{h_P}{h_P \circ\, \theta_P} - 1$
is smaller than or equal to the degree of~$h_P$~.
But according to~$(\ref{gpsi_Q'_sur_gpsi_Q})$, the zeros of~$h_P$ are also zeros of~$\frac{h_P}{h_P \circ \theta_P} - 1$. 
So we just proved that if  the function 
$\frac{h_P}{h_P \circ \theta_P} - 1$
is non-zero, it has the same divisor as  $h_P$. Therefore
there exists a  constant~$c \in k^a$ such that:
$$
\frac{h_P}{h_P \circ \theta_P} - 1 = c h_P
\qquad \text{or equivalently:} \qquad
\frac{1}{h_P \circ \theta_P} = \frac{1}{h_P} + c.
$$
Since~$\theta_P$ has finite order~$e$ and $k^a$ has characteristic zero, 
we deduce that~$ce = 0$, then~$c=0$, then~$h_P \circ \theta_P = h_P$. 
In other words, $\theta_P$ is an automorphism of the cover~$h_P = g \circ \psi \circ \chi_P : \cF_P \to \PP^1$.
According to point~(\ref{fcts_fg_Aut_gpsichi}) of lemma~\ref{fcts_fg},
we deduce that~$\theta_P$ is the  restriction to~$\cF_P$ of an element in ~$A \times B$.
We replace~$\theta$ by $\theta$ composed with  the inverse of this  element. So we now  assume that
$\theta$ acts trivially on  $\cF_P$ for some~$P$.
In particular~$\theta$ fixes every singular point on~$\cF_P$. 
So~$\theta$ stabilizes every component~$\cE_Q$. The  restriction~$\theta_Q$ of $\theta$ to~$\cE_Q$
is an automorphism of~$\chi_Q$ according to~(\ref{eq:chiq}). Further $\theta_Q$ fixes one point (and even every point)
in the unramified fiber above~$P$ of the Galois cover~$\chi_Q : \cE_Q \to \cB$. Therefore $\theta_Q$ is the identity. We have proved
that~$\theta$ is trivial on  every  component~$\cE_Q$.

To finish with, we now prove that $\theta$ is also trivial on the  components~$\cF_{P'}$
for every zeros~$P'$ of~$f$. Remind we have already assumed this to be true for one of these zeros.
We call~$\theta_{P'}$ the restriction of~$\theta$ to~$\cF_{P'}$. 
We already proved that $\theta_{P'}$ is the restriction of an element in~$A \times B$.
Further $\theta_{P'}$ fixes all the singular points
of~$\cK_0$ lying on~$\cF_{P'}$.
These points are the zeros of~$g \circ \psi \circ \chi_{P'}$. So we just need to prove that
the action of~$A \times B$ on the set of zeros of $g \circ \psi \circ \chi_{P'}$ is free.
This is certainly
the case for elements in~$B$ because the zeros of~$g \circ \psi$ are,
by hypothesis, unramified in the Galois cover~$\chi_{P'} : \cF_{P'} \to \cD$.
This is true also for elements in~$A \times B$ because the action
of~$A$ on the set of zeros of~$g \circ \psi$ is free.
\end{preuve}

\subsection{Deformations}\label{subsection:deformation}

In this paragraph we deform the two stable curves $\cJ_0$ and $\cK_0$.
If  $t\in k^a$ is a scalar, it is natural to consider the curve $\cJ_t$ drawn on the surface 
 $\cI=\cB\times \cD$ and defined by the equation
 $f(P)\times g(\psi(Q))=t$. We call $\cK_t$ the inverse image of 
 $\cJ_t$ by $\chi$.
In this paragraph and in the next one, we
 shall prove that for almost all scalars  $t$ in  $k$, the curve  $\cK_t$ is smooth, 
geometrically integral, with $k^a$-automorphism group equal to  $A\times B$, 
and with the  same field of moduli and the same fields of definition as the original
cover $\varphi$.
To this end, we would like to consider the families $(\cJ_t)_t$ and   $(\cK_t)_t$ as 
 fibrations above $\PP^1$.
We should be careful however : the family  $(\cJ_t)_t$ has base points. 
%JMC je remplace ci dessous \overline{\QQ} par  K
%trac{\'e}es sur~$(\cB \otimes_KL) \times \cD$et~$\cS$ respectivement.
So we first have to blow up $\cI=\cB\times \cD$  along 
$$
\Delta = ((f)_\infty\times (g\circ \psi)_0) \cup ((f)_0\times
(g\circ\psi)_\infty).
$$
Note that  $\Delta$ is the  union  of 
$2\times \deg(f)\times \deg(g\circ\psi)$ simple geometric points.
We denote by  $\cI_{\infty, \infty}\subset  \cI= \cB\times \cD$ 
the complementary open set of  $((f)_\infty\times \cD) \cup (\cB\times
(g\circ\psi)_\infty)$ in  $\cB\times \cD$.
We similarly define  $\cI_{0,0}$, $\cI_{0,\infty}$, $\cI_{\infty,0}$.
These four open sets cover  $\cB\times \cD$.

Let  $\PP^1=\Proj(k^a[T_0,T_1])$ be the projective line over $k^a$.
We set  $F=1/f$ and  $G=1/g$. 
Let  $\cJ_{\infty,0}\subset \cI_{\infty,0}\times 
\PP^1$ be
the set of  $(P,Q,[T_0:T_1])$ such that  $f(P)T_0=G(\psi(Q))T_1$.
Let  $\cJ_{0,\infty}\subset \cI_{0,\infty}\times 
\PP^1$ 
be the set of  $(P,Q,[T_0:T_1])$ such that $g(\psi(Q))T_0=F(P)T_1$.
Let  $\cJ_{\infty,\infty}\subset \cI_{\infty,\infty}\times 
\PP^1$ 
be the set of  $(P,Q,[T_0:T_1])$ such that  $f(P)g(\psi(Q))T_0=T_1$.
Let  $\cJ_{0,0}\subset \cI_{0,0}\times 
\PP^1$ be the set of $(P,Q,[T_0:T_1])$ such that  $T_0=F(P)G(\psi(Q))T_1$.
We glue together these four algebraic varieties and obtain
a variety  $\cJ\subset \cI\times \PP^1$. Let 
$\pi_\cI : \cJ\rightarrow \cI$ be the  projection on the first factor and let 
 $\pi_C : \cJ\rightarrow \PP^1$ be  the projection
on  $\PP^1$. This is a flat, projective, surjective morphism.

Let  $\cK\subset \cS\times \PP^1$  be the inverse image of $\cJ$ by  $\chi\times \Id$ where
 $\Id : \PP^1\rightarrow \PP^1$ is the identity.
This is the blow up of  $\cS$ along  $\chi^{-1}(\Delta)$. 
Note that  $\chi^{-1}(\Delta)$ is the  union  of 
$\deg(\chi)\times \deg(f)\times \deg(g\circ\psi)$ simple geometrical points because
$\chi$ is unramified above  $\Delta$.
Actually, $\cK$ is the  normalization of 
$\cJ$ in $k^a(S \times \PP^1)$. We  denote by  $\chi : \cK \rightarrow
\cJ$ the corresponding  morphism. 
We call  $\pi_\cS : \cK\rightarrow \cS$  the  projection on the first factor.
We call  $\pi_D : \cK \rightarrow \PP^1$
the projection on the  second factor. This is the composed morphism 
 $\pi_D=\pi_C\circ \chi$. This is a flat, proper and surjective morphism.

%PROJECTIF???????????????????????????????????????????????????????????????????????????????????????????????????????????????????????????????????????????????????????????????????????????????????????????????????

Let   $\AA^1\subset \PP^1$ be the spectrum of  $k^a[T]$
where~$T = \frac{T_1}{T_0}$.
Using the  function $T$
we identify $\PP^1(k^a)$ and $k^a \cup \{\infty \}$.
If  $t$ is a point on  $\PP^1(k^a)$  we denote  by $\cJ_t$ 
the  fiber of  $\pi_C$ above  $t$ and  $\cK_t$
the  fiber of $\pi_D$ above $t$.  The  restriction
of $\pi_\cI$ to $\cJ_t$ is a closed  immersion.
So we can see  $\cJ_t$ as a curve drawn on  $\cI=\cB\times\cD$.
Similarly, the restriction of $\pi_\cS$ to  $\cK_t$ is a closed immersion.
So we can see  $\cK_t$ as a curve drawn on  $\cS$.
In particular, the fiber of $\pi_C$ at $0$ is  isomorphic
by  $\pi_\cI$ to the stable curve  $\cJ_0$ introduced in paragraph \ref{subsection:deux}.
Similarly, the fiber of $\pi_D$ at $0$ is  isomorphic
by  $\pi_\cS$ to the stable curve  $\cK_0$ introduced in  paragraph \ref{subsection:deux}.

We call $\cJ_\eta$ the generic fiber of $\pi_C$
and  $\cK_\eta$ the generic fiber of  $\pi_D$.

We now show that the curve $\cJ_\eta$ over~$k^a(\PP^1)$
is geometrically connected  and  that for almost every $t\in \PP^1(k^a)$ the curve
$\cJ_t$ over~$k^a$ is connected.
According to Stein's factorization theorem  \cite[Chapter 5, Exercise 3.11]{Liu}, we can factor
 $\pi_C : \cJ
\rightarrow \PP^1$ as  $\pi_f\circ \pi_c$ where  $\pi_c$   has geometrically connected  fibers 
and  $\pi_f$ is finite and dominant. The fiber of  $\pi_f$ above
 $0$ is trivial because  $\cJ_0$ is  connected and reduced. Therefore
the degree of  $\pi_f$ is $1$ 
according to  \cite[Chapter 5, Exercise 1.25]{Liu}. Therefore $\pi_f$ is an  isomorphism
  above a non-empty open set of  $\PP^1$. The generic fiber  $\cJ_\eta$ is geometrically connected
over  $k^a(\PP^1)$ and for almost all  $t\in \PP^1(k^a)$ the curve  $\cJ_t$ over~$k^a$
is connected.

We now show that $\cJ_\eta$ is smooth (and therefore geometrically integral). 
Indeed, it is smooth outside the points~$(P,Q) \in \cJ_\eta\subset \cB \times \cD$
where~$df(P) =0$ and $d(g \circ \psi)(Q) = 0$. Such points are defined over $k^a$.
Therefore the function $f(P)\times g(\psi(Q))$ cannot take the transcendental value~$T$
at these points.

The ramification locus  $\Gamma\subset \cI$ of $\chi$ cuts the fiber $\cJ_0$ transversally.
Therefore it cuts the generic fiber  $\cJ_\eta$ transversally.
So  $\cK_\eta$ is smooth and geometrically integral.
Thus for almost every $t\in k^a$ the fibers  $\cJ_t$ and  $\cK_t$
are smooth and integral.

Finally, our knowledge of~$\Aut_{k^a}^{\text{adm.}} (\cK_0)$
enables us to show that~$\Aut_{{k(\PP^1)}^s}(\cK_\eta) \simeq A \times B$.
Indeed, set  $R = k^a[[T]]$ the completed local ring at the point~$T=0$ of $\PP^1$. 
%et $\eta$  le point g{\'e}n{\'e}rique de $\Spec  R$.
%le compl{\'e}t{\'e} de $\AA^1_L\otimes_Lk^a$ en $0$.
%On note $\hat \cJ = \cJ\otimes \Spec(R)$ 
%et  $\hat j : \hat \cJ \rightarrow \Spec(R)$. On note 
%$\hat \cK = \cK\otimes \Spec(R)$
%et $\hat k : \hat \cK \rightarrow \Spec(R)$. Ce sont deux courbes
%stables sur $\Spec(R)$. 
The curve  $\hat \cK = \cK \times_{\PP^1}\Spec (R)$ is stable over
the spectrum of  $R$. According to \cite[Chapter 10, Proposition 3.38, Remark 3.39]{Liu}
the functor ''automorphism group''  $t\mapsto \Aut _t(\hat \cK_t)$ 
is representable by a finite unramified scheme over $\Spec  R$ 
and the specialization morphism   $\Aut_{k^a((T))}(\hat \cK_\eta)
\rightarrow \Aut_{k^a} (\cK_0)$ is injective. According
to lemma \ref{lem:defo}, the image of this
morphism is included in the subgroup
of admissible  $k^a$-automorphisms of $\cK_0$. 
Since  $\Spec  R$ has no unramified cover, we deduce the following estimate
for the automorphism group of the generic fiber
$$
A\times B\subset \Aut_{{k^a(\PP^1)}^s}(\cK_\eta) \subset 
\Aut_{{k^a((T))}}(\hat \cK_\eta)\subset \Aut^{\text{adm.}}_{k^a}(\cK_0).
$$
We have already  seen that the rightmost group is equal to
$A\times B$. So~$\Aut_{{k^a(\PP^1)}^s}(\cK_\eta) = A \times B$ as was to be proved.

\subsection{Fields of moduli and fields of definition
of fibers}

As we have seen in paragraph~\ref{subsection:deformation},
for almost all~$t \in \AA^1(k)$, the 
fiber~$\cK_t$ is smooth and geometrically integral. 
Using lemma \ref{lem:specialisation} on the specialization
of the automorphism group in a family of curves, we deduce that
for almost all~$t \in \AA^1(k)$, the group of   $k^a$-automorphisms
of the fiber $\cK_t$ is isomorphic to
the group of  ${k(\AA^1)}^s$-automorphisms of the generic fiber.
Since the latter group is isomorphic to the automorphism
group $A\times B$  of the surface~$\cS$, we deduce
that, for almost all~$t$,  the restriction map above
is an isomorphism:
\begin{equation}\label{eq:auto}
\Aut_{k^a}(\cS) \overset{\simeq}{\longrightarrow} \Aut_{k^a}(\cK_t).
\end{equation}
Now let  $t\in k$ be such that  $\cK_t$ is smooth and geometrically
integral and such that~$\Aut_{k^a}(\cK_t)=A\times B$.
We call~$\pi_t : D_t \to S$ the corresponding
embedding.

We construct a functor~$\FF_t : \Mod_S \to \Mod_{\pi_t}$.
We first define the image of an object by~$\FF_t$.
Let~$l$ be a finite extension of~$k$
inside~$k^a$ and~$S_l$ an $l$-model of~$S$. Using the functor~$\Mod_S \to \Mod_U$ given in
section~\ref{s_def_S} followed by the functor~$\Mod_U \to \Mod_\psi$ of the proof of
lemma~\ref{lem:produitdefi}, we obtain an $l$ model~$\psi_l : Z_l \to X_l$
of the cover~$\psi$, where~$X_l = X_k \times_{\Spec(k)} \Spec(l)$ and~$Z_l$ is a $l$-model
of~$Z$. There exists also an abelian cover~$\chi_l : S_l \to X_l \times Z_l$. It is
uniquely defined up to an automorphism of~$S_l$. We denote by~$C_{t,l}$ the
curve on~$X_l \times Z_l$ defined by the equation~$f \cdot g\circ \psi_l - t = 0$. Let~$D_{t,l}$
be the inverse image of~$C_{t,l}$ by~$\chi_l$. Let~$\pi_{t,l} : D_{t,l} \hookrightarrow S_l$ be the
inclusion map. The image of the object~$S_l$ by the functor~$\FF_t$ is defined to
be~$\pi_{t,l}$. We still need to define the image of a morphism by the functor~$\FF_t$.
Let~$l'$ be another finite extension of~$k$ and let~$\sigma : l \to l'$ be a $k$-homomorphism.
Let $S'_{l'}$ be an $l'$-model of~$S$ and let~$\alpha : S_l \to S'_{l'}$ be a morphism
above~$\Spec(\sigma)$. We call~$\pi'_{t,l'} : D'_{t,l'} \hookrightarrow S'_{l'}$ the image
by~$\FF_t$ of~$S'_{l'}$. Then~$\alpha$ maps~$D_{t,l}$ to~$D'_{t,l'}$. We
denote by~$\beta$ the restriction of~$\alpha$ to~$D_{t,l}$. The image of~$\alpha$ by~$\FF_t$
is defined to be the morphism~$(\alpha,\beta)$ from~$\pi_{t,l}$ to~$\pi'_{t,l'}$.
If we compose~$\FF_t : \Mod_{S} \to \Mod{\pi_t}$ with the forgetful
functor~$\Mod_{\pi_t} \to \Mod_{D_t}$, we obtain a cartesian
functor~$\GG_t : \Mod_S \to \Mod_{D_t}$.
Further, identity~(\ref{eq:auto}) implies that the functor~$\GG_t$
is fully faithful. Therefore, by proposition~\ref{obstruction_transportor},
both~$S$ and~$\cK_t$ have $k$ as field of moduli and a $k$-extension is a field of
definition of~$S$ if and only if it is a field of definition of~$\cK_t$.
In view of section~\ref{s_def_S}, $\cK_t$,~$\psi$ and~$\varphi$ also share the same
fields of definition. Theorem~\ref{theoreme:courbe} is proved.

\section{Six lemmas about curves and surfaces}

In this section we state and prove seven lemmas that are needed in the
proof of theorem \ref{theoreme:courbe}.

\subsection{About curves and products of two curves}

\begin{lem}\label{lem:adjonction}
Let  $k$ be a algebraically closed field.
Let   $\cB$ and  $\cC$  be two  projective, smooth and integral curves over $k$.
Let $\beta$ be the genus of $\cB$
and let $\gamma$ be the genus of $\cC$.
We fix a geometric point  $P$ on  $\cB$ and a geometric point  $Q$ on  $\cC$.
We identify the curves $\cB$ and $\cB\times Q$ and the curves~$\cC$ and $P\times \cC$.
Let   $\Gamma$ be a  divisor on  $\cB\times \cC$ of bidegree~$(b,c)$,
i.e. $b=\cB\cdot\Gamma$ and
 $c=\cC\cdot\Gamma$. The  virtual arithmetic genus
  $\pi$ of  $\Gamma$ is at most  $1+bc+c(\beta-1)+b(\gamma-1)$. When
$b=c$ this bound reads  $1+2b(\beta -1)+b^2$.
\end{lem}

\begin{preuve}
We follow the lines of  Weil's proof of the  Riemann hypothesis 
for curves (cf.~\cite[Exercise~V-1.10]{Hartshorne}).

The algebraic equivalence class of the canonical divisor
on $\cB\times\cC$ is  $K=2(\beta -1)\cC+2(\gamma -1)\cB$. 
We recall that the virtual arithmetic genus~$\pi$, as defined in~\cite[Exercise~V-1.3]{Hartshorne},
is such that~$\pi=\frac{D\cdot(D+K)}{2}+1$ for every divisor~$D$.
Thus:
$$
\pi = \frac{D\cdot\left(D + 2(\beta -1)\cC+2(\gamma -1)\cB\right)}{2} + 1
    = \frac{D \cdot D + 2c(\beta - 1) + 2b(\gamma - 1)}{2} + 1,
$$
and we just need to bound the self intersection~$D \cdot D$. 
We deduce from Castelnuovo's and Severi's inequality
 (cf.\cite[Exercise~V-1.9]{Hartshorne}) that~$D \cdot D \leq 2bc$. This
finishes the proof of the lemma.
\end{preuve}

\begin{lem} \label{lem:fibrescoupenttransversalement}
Let  $k$ be an algebraically closed field.
Let~$X$ and~$Y$ be two  projective, smooth,
integral curves over $k$.
Let  $\Gamma$ be an effective divisor without multiplicity on the 
surface~$X \times Y$. 
Let  $f \in k(X)$ be a non-constant   function.
For all but finitely many scalars $x$ in $k$, the
divisor~$(f)_x \times Y$ crosses  transversally~$\Gamma$, 
where~$(f)_x$ is the positive part of the  divisor of $f-x$.
\end{lem}

\begin{preuve}
We call~$p_\cB : \cB \times \cC \to \cB$ the projection
on the first factor. Let~$E$ be the set of  points
in~$\cB(k)$ such that at least one of the following condition
holds: $p_\cB^{-1}(P)$ contains a singular  point on~$\Gamma$, or
$p_\cB^{-1}(P)$ contains
a ramified point of the  morphism~$p_\cB : \Gamma \to \cB$, 
or the  fiber~$p_\cB^{-1}(P)$ is contained in~$\Gamma$.
The set $E$ is finite. For all  $x\in k$ but
finitely many, the fiber  $f^{-1}(x)$ avoids $E$ and it 
is  simple. 
\end{preuve}

\begin{lem}\label{lem:autofibres}
Let  $k$ be an algebraically closed field.
Let  $\cB$ be a projective, smooth, integral curve
over   $k$. Assume the genus of  $\cB$ is at least  $2$.
Let  $f\in k(X)$ be a  non-constant 
function. We note  $G$ the  group of $k$-automorphisms of $f$.
This is the set of all  $k$-automorphisms $\theta$ 
of $\cB$  such that  $f \circ \theta = f$. For any
 $x\in \PP^1 (k)$,  we  note 
$(f)_x=f^{-1}(x)$ the fiber above  $x$
and $G_x$ the group of $k$-automorphisms
of  $\cB$ that stabilize the set of 
$k$-points of $(f)_x$.

For all  $x$ in $\PP^1 (k)$ but finitely many we have  $G_x=G$.
\end{lem}

\begin{preuve}
The group  $H=\Aut_{k}(\cB)$ of $k$-automorphisms
of $\cB$ is finite because the genus of $\cB$ is at least two.
Let  $\theta$ be an  automorphism in $H \setminus G$ and let 
  $x\in \PP^1 (k)$. Assume that the  $k$-points in
 $(f)_x$ are permuted by  $\theta$. Let
 $P$ be one of them. Then $f\circ \theta (P)=f(P)=x$.
So  $P$ is a zero of the
non-zero function  $f\circ \theta -f$.
For each  $\theta$ there are finitely many such
zeros. And the  $\theta$ are finitely many.
 So the images by $f$ of such $P$'s are finitely many also.
\end{preuve}

\begin{lem}\label{lem:fonctionsauto}
Let  $k$ be a field. 
Let  $\cB_k$ be a projective, smooth, geometrically integral curve
over $k$. Set~$\cB = \cB_k \times_{\Spec(k)} \Spec(k^a)$ and assume
that $\cB$ has genus at least  $2$.
Let  $\cC$ be a projective, smooth, integral
curve over $k^a$ and let  $\varphi : \cC\rightarrow \cB$ be a  non-constant 
 $k^a$-cover.
If~$f$ is any non-constant function in $k^a(\cB)$
then $\Aut(\varphi) \subset \Aut(f \circ \varphi)$.
Let  $V\subset k(\cB_k)$ be the set of functions
$f \in k(\cB_k)$  such that $\Aut(\varphi) \not= \Aut(f \circ \varphi)$.
This set  $V$ is contained in a finite
 union of strict $k$-subalgebras
of  $k(\cB_k)$.
\end{lem}

\begin{preuve}
The statement to be proven concerns the three function fields
~$k^a(f) \subset k^a(\cB) \subset k^a(\cC)$ and the groups involved
are the following ones:
$$
\begin{cases}
\Aut(\varphi) = \Aut_{k^a(\cB)}(k^a(\cC)), \\
\Aut(f \circ \varphi) = \Aut_{k^a(f)}(k^a(\cC)), \\
\Aut(\cC) = \Aut_{k^a}(k^a(\cC)),
\end{cases}
\quad\Rightarrow\quad
\Aut(\varphi) \subset \Aut(f\circ\varphi) \subset \Aut(\cC).
$$
Now, the set~$V$ can be described as follows:
$$
V
=
\left(
\bigcup_{\theta \in \Aut(Y) \setminus \Aut(\varphi)}
k^a(\cC)^\theta \cap k^a(\cB)
\right) \cap k(\cB_k)
=
\bigcup_{\theta \in \Aut(\cC) \setminus \Aut(\varphi)}
k^a(\cC)^\theta \cap k(\cB_k).
$$
This is a  union of sets indexed by elements in the finite
set~$\Aut(\cC)\setminus \Aut(\varphi)$ (remind $\Aut(\cC)$ is finite because the genus of~$\cC$ 
is at least $2$). Since~$\theta \not \in \Aut(\varphi)$, 
each~$k^a(\cC)^\theta \cap k^a(\cB)$ 
is a strict subfield  of~$k^a(\cB)$ containing $k^a$. 
Therefore~$k^a(\cC)^\theta \cap k(\cB_k) \subsetneq k(\cB_k)$.
\end{preuve}

\subsection{Deformation of an automorphism of a nodal curve}

In this paragraph we give
a {\it necessary}
condition for extending an automorphism of a nodal curve to 
a given deformation of this curve.

Let~$R$ be a complete 
discrete valuation ring. Let~$\pi$ be a uniformizing parameter
and let~$k$ be the residue field. We assume $k$ is  algebraically closed.
 Let~$\cK$ be a semi-stable curve over~$\Spec(R)$. We note
$\cK_\eta$ the generic fiber and $\cK_0$ the special fiber.
We assume $\cK_\eta$ is smooth over the fraction field of $R$.
Let~$T$ be a singular point of~$\cK_0$. According 
to~\cite[Chapter~10, Corollary~3.22]{Liu}, the completion of the local ring
of~$\cK$ at~$T$ takes the form:
$$
\hat{\cO}_{\cK, T} = R[[f,g]]/\langle f g - \pi^e \rangle
$$
where~$e$ is a positive integer. 
This integer is called the  {\em thickness} of~$\cK$ at~$T$. We also say that~$f$ and~$g$ form
a coordinate system for~$\cK$ at~$T$. If we reduce modulo $\pi$, we obtain the completion
of the local ring of~$\cK_0$ at~$T$:
$$
\hat{\cO}_{\cK_0, T} = \hat{\cO}_{\cK, T} / \langle \pi \rangle
= k[[\overline{f},\overline{g}]]/\langle \overline{f} \overline{g} \rangle,
$$
where~$\overline{f} = f \bmod{\pi}$ and~$\overline{g} = g \bmod{\pi}$.
Because $T$ is an ordinary double point,~$\cK_0$ has two branches $\cF$ and $\cG$ at~$T$.
These correspond to
  the two irreducible components of the completion at $T$. Be careful  that these two branches may
lie on the same irreducible component of $\cK_0$. 
Anyway, the functions~$\overline{f}$ and~$\overline{g}$ are the uniformizing parameters
of either branches. We call~$P$ and~$Q$ the  points of~$\cF$
and~$\cG$ above~$T$.

Now let~$T'$ be another singular point of~$\cK_0$, and let~$f'$, $g'$, $e'$, $\cF'$, and $\cG'$ 
the corresponding data.

Let~$\theta$ be an  automorphism of~$\cK$ over~$R$ such that~$\theta(T) = T'$ 
and~$\theta(\cF) = \cF'$, $\theta(\cG) = \cG'$. One easily checks that
the functions~$f' \circ \theta$ and~$g' \circ \theta$ form a coordinate
system for~$\cK$ at~$T$. We deduce that~$e' = e$ and that both~$f' \circ \theta / f$ and~$g' \circ \theta / g$
are units in~$\hat{\cO}_{\cK, T}$ (indeed, in either fraction,  the numerator and denominator have
the same Weil divisor). Since~$f \times g = \pi^e = f' \circ \theta \times g' \circ \theta$, we 
have~$\frac{f' \circ\, \theta}{f}(T) \times \frac{g' \circ\, \theta}{g}(T) = 1$.
We reduce this identity modulo~$\pi$ and obtain the following identity
where the first factor is a function on $\cF$ evaluated at $P$ while the second factor
is a function on $\cG$ evaluated at $Q$:

\begin{equation}\label{eq:admissiblegene}
\frac{\overline{f}' \circ \overline{\theta}}{\overline{f}}(P) \times \frac{\overline{g}' \circ \overline{\theta}}{\overline{g}}(Q) = 1.
\end{equation}

This leads us to the following definition.

\begin{definition}\label{defi:deform}
Let  $R$ be a complete discrete valuation ring. Assume that  the residue field 
$k$ is algebraically closed. Let $\cK$ be a semi-stable curve over $\Spec(R)$. 
The generic fiber of  $\cK$ is assumed to be smooth. Assume we are given
a coordinate system at each singular point of the special fiber $\cK_0$.
Let  $\btheta$ be an   automorphism of the special fiber $\cK_0$. 
We say that  $\btheta$ is {\it admissible} in $\cK/\Spec(R)$
if for every singular point $T$ of $\cK_0$, the image  $\btheta(T)$
has the same thickness as  $T$ in $\cK$, and if equality  (\ref{eq:admissiblegene})
holds true.
\end{definition}

We have just proved the following lemma.

\begin{lem}\label{lem:defo}
With the notation of definition \ref{defi:deform}, the set of  automorphisms 
of $\cK_0$  that are admissible in  $\cK/\Spec(R)$  
form a subgroup of  $\Aut_k(\cK_0)$.
If  $\theta$ is an automorphism of $\cK$ over $\Spec(R)$, its reduction
$\btheta = \theta \bmod \pi$ is an automorphism of $\cK_0$
and  is admissible in  $\cK/\Spec(R)$.
\end{lem}

One may compare this statement with  \cite[Theorem 3.1.1]{Wewers}
where the  deformation of morphisms between two
distinct curves is studied.

\begin{remarque}
It must be pointed out that the converse of lemma \ref{lem:defo} is not true.
For example, consider
the  elliptic curve  $E$ with modular invariant $j=0$ (or $1728$). Every automorphism
of $E$
is admissible because there are no singular  points on the curve (the condition in definition
\ref{defi:deform} is empty). However, the only
automorphisms that can be extended to the generic elliptic curve are the identity and the involution.
\end{remarque}

\subsection{Automorphisms of curves in a family}

%ici j'ai travaill{\'e} un peu
In this section we state and prove a lemma
of specialization of the automorphism group
of  curves in a family.

\begin{lem}\label{lem:specialisation}
Let  $k$ be a field  and let  $U$ be
a smooth, geometrically integral curve over $k$. Let
$X$ be a quasi-projective, smooth, geometrically integral
surface over $k$. Let 
$\pi : X \rightarrow U $ be a surjective, projective,
smooth  morphism of relative dimension  $1$. Assume
that for any point   $x$
of~$U$, the  fiber $X_x$ at $x$
  is geometrically integral.
We call  $\eta$ the generic point of  $U$ and 
$\bar X_\eta = X_\eta\times_{\Spec(k(U))}\Spec({k(U)}^a)$ 
the generic fiber, seen as a curve
over the algebraic closure of the function field of the
basis $U$. We assume the genus of $X_\eta$ is at least $2$.

There exists a non-empty  open subset $V$ of $U$ over $k$
such that for any geometric  point $x \in V(k^a)$ the group
of $k^a$-automorphisms of the fiber at $x$ is equal 
to the group 
$\Aut_{{k(U)}^a}(\bar X_\eta )$ of  automorphisms of  $\bar X_\eta$.
\end{lem}

The following proof was communicated to us by Qing Liu.

\begin{preuve}
This is a consequence of a
general result by Deligne-Mumford. 
Let $X \rightarrow  S$ be a flat projective
morphism  over  a noetherian scheme $S$. The
functor $T \rightarrow   \Aut_T(X_T)$
from the category of  $S$-schemes to the category
of groups is representable
by a group scheme
$\Aut_{X/S}$ over $S$. See \cite[Exercise 1.10.2]{kol} for example.

When $X\rightarrow S$ is a stable curve with genus at least
$2$,  Deligne and Mumford \cite[Thm 1.11]{DM} prove that the scheme
$\Aut_{X/S}$ is finite and unramified over $S$.

In our lemma,  $S$ is a the smooth, geometrically integral curve $U$ over $k$. Replacing
$S$ by  a non-empty
open subset,
we may assume that
$\Aut_{X/S}$ is finite étale over $S$.
At the expense of   a finite
surjective
base change $T\rightarrow  S$,
we may assume that the generic fiber
of $\Aut_{X/S}\rightarrow S$
consists of rational points. So $\Aut_{X/S} \rightarrow S$
is now a disjoint union of étale sections and the fibers
have constant degree.
In particular, the fibers are constant
and the specialization maps
$\Aut_S(X)=\Aut_{X/S}(S)\rightarrow \Aut_s(X_s)=\Aut_{X/S}(k(s))$
are isomorphisms.
\end{preuve}

\bibliographystyle{alpha}
%\bibliography{moduli}

\end{document}